\documentclass{amsart}

\usepackage{amssymb, amscd, latexsym, graphicx, psfrag}
\usepackage[all]{xy}

\newtheorem{thm}{thm}
\newtheorem{theo}[thm]{Theorem}
\newtheorem{cor}{cor}
\newtheorem{coro}[cor]{Corollary}

\newtheorem{dummy}{dummy}[section]
\newtheorem{lemma}[dummy]{Lemma}
\newtheorem{theorem}[dummy]{Theorem}

\newtheorem{corollary}[dummy]{Corollary}
\newtheorem{proposition}[dummy]{Proposition}
\theoremstyle{definition}
\newtheorem{definition}[dummy]{Definition}
\newtheorem{example}[dummy]{Example}
\newtheorem{remark}[dummy]{Remark}
\newtheorem{remarks}[dummy]{Remarks}

\newcommand{\bC}{\mathbb{C}}

\newcommand{\bL}{\mathbb{L}}
\newcommand{\bP}{\mathbb{P}}
\newcommand{\bQ}{\mathbb{Q}}
\newcommand{\bR}{\mathbb{R}}
\newcommand{\bZ}{\mathbb{Z}}

\newcommand{\cA}{\mathcal{A}}
\newcommand{\cC}{\mathcal{C}}
\newcommand{\cD}{\mathcal{D}}

\newcommand{\cF}{\mathcal{F}}

\newcommand{\cK}{\mathcal{K}}
\newcommand{\cL}{\mathcal{L}}
\newcommand{\cP}{\mathcal{P}}
\newcommand{\cQ}{\mathcal{Q}}
\newcommand{\cO}{\mathcal{O}}

\newcommand{\cT}{\mathcal{T}}
\newcommand{\cU}{\mathcal{U}}
\newcommand{\cY}{\mathcal{Y}}
\newcommand{\cR}{\mathcal {R}}

\newcommand{\Hom}{\mathrm{Hom}}

\newcommand{\Pic}{\mathrm{Pic}}

\newcommand{\si}{\sigma}
\newcommand{\Si}{\Sigma}

\newcommand{\vc}{ {\vec{c}} }

\newcommand{\dbar}{\bar{\partial}}
\newcommand{\define}{\stackrel{\mathrm{def}}{=} }

\newcommand{\pa}{\partial}
\newcommand{\tri}{\triangle}
\newcommand{\LS}{ {\Lambda_\Sigma} }

\newcommand{\Spec}{\mathrm{Spec}}

\newcommand{\Sh}{\mathit{Sh}}
\newcommand{\naive}{\mathit{naive}}

\newcommand{\Perf}{\cP\mathrm{erf}}

\renewcommand{\SS}{\mathit{SS}}

\begin{document}

\title{T-Duality and Homological Mirror Symmetry of Toric Varieties}

\author{Bohan Fang}
\address{Bohan Fang, Department of Mathematics, Northwestern University,
2033 Sheridan Road, Evanston, IL  60208}
\email{b-fang@math.northwestern.edu}

\author{Chiu-Chu Melissa Liu}
\address{Chiu-Chu Melissa Liu, Department of Mathematics, Columbia University,
2990 Broadway, New York, NY 10027}
\email{ccliu@math.columbia.edu}

\author{David Treumann}
\address{David Treumann, Department of Mathematics, Northwestern University,
2033 Sheridan Road, Evanston, IL 60208}
\email{treumann@math.umn.edu}

\author{Eric Zaslow}
\address{Eric Zaslow, Department of Mathematics, Northwestern University,
2033 Sheridan Road, Evanston, IL  60208}
\email{zaslow@math.northwestern.edu}

\begin{abstract}
Let $X_\Sigma$ be a complete toric variety.
The coherent-constructible correspondence $\kappa$ of \cite{FLTZ} equates
$\Perf_T(X_\Sigma)$ with a subcategory
$Sh_{cc}(M_\bR;\LS)$ of constructible sheaves on a vector space $M_\bR.$
 The microlocalization equivalence $\mu$ of
\cite{NZ,N} relates these sheaves to a subcategory $Fuk(T^*M_\bR;\LS)$
of the Fukaya category of the cotangent $T^*M_\bR$. When $X_\Si$ is nonsingular,
taking the derived category yields an equivariant version of homological
mirror symmetry, $DCoh_T(X_\Si)\cong DFuk(T^*M_\bR;\LS)$, which is an equivalence
of triangulated tensor categories. 

The nonequivariant coherent-constructible correspondence 
$\bar{\kappa}$ of \cite{T} embeds $\Perf(X_\Si)$ into a subcategory 
$Sh_c(T_\bR^\vee;\bar{\Lambda}_\Si)$  of constructible sheaves on a compact torus $T_\bR^\vee$.
When $X_\Si$ is nonsingular, the composition of $\bar{\kappa}$ and microlocalization
yields a version of homological mirror symmetry,
$DCoh(X_\Sigma)\hookrightarrow DFuk(T^*T_\bR;\bar{\Lambda}_\Si)$,
which is a full embedding of triangulated tensor categories.

When $X_\Si$ is nonsingular and projective, the composition $\tau=\mu\circ \kappa$ is
compatible with T-duality, in the following
sense.  An equivariant ample line bundle $\cL$ has a hermitian metric invariant under the
real torus, whose connection defines a family of flat line bundles over the real torus orbits.
This data produces a T-dual Lagrangian brane $\mathbb L$ on the universal cover $T^*M_\bR$ of
the dual real torus fibration.  We prove $\mathbb L\cong \tau(\cL)$ in $Fuk(T^*M_\bR;\LS).$
Thus, equivariant homological mirror symmetry is determined by T-duality.
\end{abstract}

\maketitle

{\small \tableofcontents}
\section{Introduction}

In this paper we derive equivariant and nonequivariant versions of the
homological mirror symmetry for nonsingular complete toric varieties
from the coherent-constructible correspondence  
\cite{FLTZ, T} and microlocalization \cite{NZ,N}. The composition of 
the coherent-constructible correspondence and microlocalization sends an equivariant (resp. nonequivariant)
coherent sheaf on the toric variety to an object in the Fukaya category of the cotangent bundle of 
a vector space (resp. a compact torus).

For nonsingular projective toric varieties, the equivariant homological
mirror symmetry is determined by equivariant ample line bundles.
We prove that the image of an equivariant ample line bundle agrees
up to isomorphism with the Lagrangian constructed by T-duality.

\subsection{Main Results} \label{sec:main-results}

In this paper, we work over $\bC$. Let $X_\Sigma$ be an $n$-dimensional complete toric variety defined by
a finite complete fan $\Sigma\subset N_\bR$.
Then $T\cong (\bC^*)^n$ acts on $X_\Sigma$, and $N_\bR\cong \bR^n$ can
be identified with the Lie algebra of the maximal compact
subgroup $T_\bR \cong U(1)^n$. Let $M_\bR = \Hom_\bR(N_\bR,\bR) \cong \bR^n$
be the dual real vector space of $N_\bR$. Then
the lattice $M=\Hom(T_\bR, U(1))\cong \bZ^n$ naturally sits in $M_\bR$, and the quotient
$T_\bR^\vee= M_\bR/M$ is the dual torus of $T_\bR$.

\bigskip

The first main theorem in this paper is the following:
\begin{theo}\label{main-I}
Let $X_\Si$ be a complete toric variety defined by a finite complete fan
$\Si\subset N_\bR$.
Then there is a quasi-equivalence of $A_\infty$ categories:
\begin{equation}\label{eqn:main-I}
\tau: \Perf_T(X_\Si)\stackrel{\cong}{\longrightarrow} Fuk(T^*M_\bR;\LS).
\end{equation}
This functor intertwines the usual monoidal product on $\Perf_T(X_\Si)$ and a product structure $\diamond$ on $Fuk(T^*M_\bR;\LS)$ up to a quasi-isomorphism.
\end{theo}
In \eqref{eqn:main-I}, $\Perf_T(X_\Si)$ is the dg category of equivariant perfect
complexes on $X_\Si$ (see Section \ref{sec:schemes} for the precise definition),
and $Fuk(T^*M_\bR;\LS)$ is a subcategory  of the unwrapped Fukaya category $Fuk(T^*M_\bR)$ determined by
the fan $\Si$ (see Section \ref{sec:HMS} for the precise definition).
As we will explain in Section \ref{sec:HMS}, Theorem \ref{main-I} follows
from the results in \cite{FLTZ, NZ, N}. We use the results in \cite{N} to define
the product $\diamond$ on $Fuk(T^*M_\bR)$ (actually any cotangent bundle of a Lie group); for cotangent fibers we have
$T_{x_1} M_\bR \diamond T_{x_2} M_\bR = T_{x_1+x_2}M_\bR$. The monoidal product structure on $\Perf_T(X_\Si)$ comes from the usual tensor product of vector bundles.

When $X$ is nonsingular, taking $H^0$ of \eqref{eqn:main-I} yields the following:
\begin{coro}[equivariant homological mirror symmetry of toric varieties]\label{equivariant}
Let $X_\Si$ be a nonsingular complete toric variety   defined by a nonsingular finite complete fan
$\Si\subset N_\bR$. (In particular, $X_\Si$ is a  compact complex manifold.)
Then there is an equivalence of tensor triangulated categories:
\begin{equation}\label{eqn:equivariant}
H(\tau): DCoh_T(X_\Si)\stackrel{\cong}{\longrightarrow} DFuk(T^*M_\bR;\LS).
\end{equation}
\end{coro}
In \eqref{eqn:equivariant}, $DCoh_T(X_\Si)$ is the bounded derived category of equivariant
coherent sheaves on $X_\Si$. The equivalence \eqref{eqn:equivariant} preserves
the tensor product, so it is a stronger equivalence
than the equivalence in the usual homological mirror symmetry. Note that we do not assume
$X_\Si$ is projective in Corollary \ref{equivariant}, so a priori the other
direction of homological mirror symmetry (involving the Fukaya category of the toric
variety) does not make sense.

\bigskip

Our second main theorem concerns the nonequivariant version of Theorem \ref{main-I}:
\begin{theo}\label{main-II}
Let $X_\Si$ be a complete toric variety defined by a finite complete fan
$\Si\subset N_\bR$.
Then there is a quasi-embedding of $A_\infty$ categories:
\begin{equation}\label{eqn:main-II}
\bar{\tau}: \Perf(X_\Si) \longrightarrow Fuk(T^*T_\bR^\vee;\bar{\Lambda}_\Si).
\end{equation}
The functor $\bar \tau$ intertwines the product $\diamond$ on $Fuk(T^*T_\bR^\vee;\bar{\Lambda}_\Si)$ and the usual monoidal product on $\Perf(X_\Si)$.
\end{theo}
In \eqref{eqn:main-II}, $\Perf(X_\Si)$ is the dg category of perfect
complexes on $X_\Si$ (see Section \ref{sec:schemes} for the precise definition),
and $Fuk(T^*T_\bR^\vee;\bar{\Lambda}_\Si)$ is a subcategory  of the unwrapped Fukaya category 
$Fuk(T^*T_\bR^\vee)$ determined by the fan $\Si$ (see Section \ref{sec:HMS} for the precise definition). The diamond product $\diamond$ is similarly defined as in the equivariant case from the Lie group structure on $T^\vee_\bR$.
As we will explain in Section \ref{sec:HMS}, Theorem \ref{main-II} follows
from the results in \cite{T, NZ, N}. We conjecture that \eqref{eqn:main-II} is
a quasi-equivalence.

When $X$ is nonsingular, taking $H^0$ of \eqref{eqn:main-II} yields the following:
\begin{coro}[homological mirror symmetry for toric varieties]
Let $X_\Si$ be a nonsingular complete toric variety defined by a nonsingular finite complete fan $\Si\subset N_\bR$.
Then there is a full embedding of tensor triangulated categories:
\begin{equation}\label{eqn:nonequivariant}
H^0(\bar{\tau}): DCoh(X_\Si)\longrightarrow DFuk(T^*T_\bR^\vee ;\bar{\Lambda}_\Si)
\end{equation}
\end{coro}
In \eqref{eqn:nonequivariant}, $DCoh(X_\Si)$ is the bounded derived category of
coherent sheaves on $X_\Si$. We conjecture that \eqref{eqn:nonequivariant} is an equivalence.
In Appendix \ref{sec:appothers}, we will comment  on the relationships
among the Fukaya categories in Theorem \ref{main-I},  the physical/traditional
mirror of toric  varieties (Landau-Ginzburg/Fukaya-Seidel category),
and the relative Fukaya category in Abouzaid's work \cite{Ab1, Ab2}, when $X_\Si$ is
a nonsingular projective toric variety.

\bigskip

Our third main theorem relates Theorem \ref{main-I} to T-duality.
When $X_\Sigma$ is nonsingular and projective, we perform
an equivariant version of T-duality: for any equivariant
line bundle $\cL_{\vc}$ with a $T_\bR$-invariant hermitian metric $h$,
we construct a Lagrangian $\bL_{\vc,h}\subset T^*M_\bR$, which projects to a Lagrangian
$\bar{\bL}_{\vc,h} \subset T^*T_\bR^\vee$. We prove the following:
\begin{theo}[equivariant homological mirror symmetry is T-duality] \label{main-III}
Let $X_\Sigma$ be a nonsingular projective toric variety defined by
a fan $\Si\subset N_\bR$. For any
equivariant ample line bundle $\cL_\vc$ with an admissible\footnote{A hermitian
metric $h$ on an ample line bundle is {\em admissible} if it is real analytic,
$T_\bR$-invariant, and defines a unitary connection whose curvature is a
nondegenerate closed 2-form.} hermitian metric $h$, the T-dual Lagrangian $\bL_{\vc,h}$
(constructed in Section \ref{sec:tduality}) is an object in $Fuk(T^*M_\bR;\LS)$ and
$$
\bL_{\vc,h}\cong \tau(\cL_\vc),
$$
where $\tau$ is as in Theorem \ref{main-I}. 
\end{theo}
\noindent
By Theorem \ref{ampgen}, when $X_\Sigma$ is nonsingular and projective,
$\Perf_T(X_\Si)$ is generated by equivariant ample line bundles.
Therefore equivariant homological  mirror symmetry
\eqref{eqn:equivariant} is determined by T-duality. Theorem \ref{main-I}
and Theorem \ref{main-III} imply the following.
\begin{coro}[subcategory generated by T-dual Lagrangians]
Let $X_\Si$ be a nonsingular projective toric variety defined
by a fan $\Si\subset M_\bR$. Then $Fuk(T^*M_\bR;\LS)$ is generated by
the T-dual Lagrangians $\bL_{\vc,h}$ of equivariant ample line
bundles $\cL_\vc$ on $X_\Si$.
\end{coro}

\subsection{Simple Example}
The simple example of $\bP^1=\bC\cup\{\infty\}$ is instructive.
The $\bC^*$ action is $t: z \mapsto t\cdot z.$  Write $z= e^{y+\sqrt{-1}\theta},$ so
$\theta\in S^1$ coordinatizes the real torus orbit.
The divisors $p_0=0$ and $p_\infty= \infty$ span the equivariant
Picard group.   The equivariant line bundle
$\cO_{\bP^1}(a p_0 +b p_\infty),$ $a,b\in \bZ,$ admits an $S^1$-invariant
hermitian metric $h = \frac{|z|^{2b}}{(1+|z|^2)^{a+b}}$ and
associated connection 1-form $A = \frac{1}{\sqrt{-1}}\partial_y \log h \,d\theta$.
On each real torus $y = \text{\em const},$ this connection has monodromy determined by the value of $\gamma = -\partial_y \log h\vert_y,$
a coordinate on the dual $S^1$.  Letting $y$ vary determines a submanifold $\mathbb L = \{(y,\gamma) \;|\; \gamma = -\partial_y \log h\}\subset \bR^2.$
By the explicit form of $h,$ we find
$\gamma = \frac{(a+b)e^{2y}}{1+e^{2y}} - b.$  The
{\sl non}equivariant
bundle is $\cO_{\bP^1}(a+b),$ and note that keeping the sum $a+b$
fixed and varying $b$ amounts to lattice translations in the
universal cover $\bR$ of the dual torus $S^1.$  Inverting equations,
we can write ${\mathbb L}$ as a graph over an interval over length $|a+b|$ in $\bR,$ which corresponds to a constructible sheaf by \cite{NZ}.

\smallskip
\begin{center}
\includegraphics[scale=0.8]{fig.1}
\parbox{10cm}{Fig.1 Our procedure for $\mathbb
P^1$. The three Lagrangians shown above come from equivariant line
bundles $\mathcal O_{\bP^1}(p_0)$, $\mathcal O_{\bP^1}(p_\infty)$ and $\mathcal O_{\bP^1}(-2p_0)$.
By microlocalization \cite{NZ}, they correspond (up to shifts) to three constructible sheaves on $\bR$:
$i_{(0,1)!} \bC_{(0,1)}$, $i_{(-1,0)!}\bC_{(-1,0)},$ and $i_{(-2,0)*}\bC_{(-2,0)}$, respectively, where
$i$ is the inclusion of the indicated open interval into $M_\bR\cong\bR$.
}
\end{center}

\subsection{Relation to the work of others}
\label{others}

The present work is much related to results of several authors.  Below are some comparisons;
further details are given in Appendix \ref{sec:appothers}.

Homological mirror symmetry for toric Fano varieties was conjectured by Kontsevich \cite{K2}.
A physical proof of mirror symmetry was given by Hori-Vafa \cite{HV}.
The mirror of a toric Fano manifold is a Landau-Ginzburg model $((\bC^*)^n, W)$ where the superpotential $W:(\bC^*)^n\to \bC$ is a holomorphic function.
The homological mirror conjecture states (in one direction)
that the derived category of coherent sheaves on the toric Fano manifold is equivalent to the derived Fukaya-Seidel category $FS((\bC^*)^n, W)$ of the Picard-Lefschetz fibration defined by $W.$

Seidel proves homological mirror symmetry for $\bP^2$ in \cite{S2}.
Auroux-Katzarkov-Orlov prove it for weighted projective planes and their
noncommutative deformations in \cite{AKO1}, and
for (not necessarily toric) del Pezzo surfaces in \cite{AKO2}.
Ueda proves it for toric del Pezzo surfaces \cite{U};
Ueda-Yamazaki prove it for toric orbifolds of toric del Pezzo surfaces.
Bondal and Ruan \cite{BR} announced a proof of homological
mirror symmetry for weighted projective spaces,
generalizing the result by Auroux-Katzarkov-Orlov
on weighted projective planes \cite{AKO1}.

The version here is somewhat different, but conjecturally related (see Section \ref{sec:abouzaid})
and much closer to Abouzaid's work \cite{Ab1, Ab2}.
Torus equivariance is encoded in the $\bZ^n$ grading  of morphisms in various categories introduced in \cite{Ab2}.

Recently, Subotic constructed a monoidal structure on the extended Fukaya category of
any Lagrangian torus fibration with a section \cite{Su}.
%
%

We also mention some complementary work which studies the A-model for
toric varieties \cite{CO,FOOO1,FOOO2} and varieties with effective anticanonical divisors \cite{Au}, and
the relation to the Landau-Ginzburg mirror---especially Chan-Leung
\cite{CL}, which employs similar T-duality reasoning.

\subsection{Outline}
Section \ref{sec:notation} contains notation and conventions for categories, sheaves and toric varieties.
In Section \ref{sec:HMS}, we derive Theorem \ref{main-I} and Theorem \ref{main-II}.
In Section \ref{sec:tduality}, we perform an equivariant version of the T-duality, and relate the resulting T-dual Lagrangians to classical objects in symplectic geometry.
%
%
In Section \ref{sec:fukaya}, we prove Theorem \ref{main-III}.
Appendix \ref{sec:appCset} contains a brief review of analytic-geometric categories
and a proof of Proposition \ref{thm:cset}.
We show that $\Perf_T(X)$ is generated by
$T$-equivariant ample line bundles in Appendix \ref{sec:ample-generate}.
We discuss the relation to the work of others in Appendix \ref{sec:appothers}.

%



\subsection*{Acknowledgments}
We thank M. Abouzaid, A. Bondal and P. Seidel for explaining relevant aspects of their work,
and D. Nadler for helpful suggestions.
The work of EZ is supported in part by NSF/DMS-0707064.
BF and EZ would like to thank the Kavli Institute for Theoretical
Physics and the Pacific Institute for the Mathematical Sciences,
where some of this work was performed.

\section{Notation and Convention}
\label{sec:notation}

\subsection{Categories}\label{sec:categories}


Throughout, we consider dg and more generally $A_\infty$ categories.
Unless otherwise stated, a category $C$ is assumed to be completed to its
triangulated envelope.  (Recall that dg and $A_\infty$ categories have
canonical triangulated structures\footnote{A triangle in an $A_\infty$ category
$C$ is distinguished if it induces a distinguished triangle in the cohomology category $H(C).$}
and completions\footnote{Here is a construction of the unique-up-to-isomorphism
triangulated envelope. The Yoneda embedding
$\mathcal Y:C\rightarrow mod(C)$
maps an object $L$ of a category $C$ to the  $A_\infty$ right $C$-module $hom_{C}(-, L)$.
The  functor $\mathcal Y$ is a quasi-embedding of $C$ into the triangulated category $mod(C)$.
Then the triangulated completion
$Tr(C)$ is the category of twisted complexes of representable modules in $mod(C)$.} \cite{S3}.)
%

%

\subsection{Schemes and coherent sheaves}
\label{sec:schemes}

All schemes that appear will be over $\bC$. If $X$ is a scheme, then we let $\cQ^\naive$  denote the dg category of
bounded complexes of quasicoherent sheaves on $X$, and we let $\cQ(X)$
denote the localization of this category  with respect to acyclic complexes
(see \cite{Dr} for localizations of dg categories).  If $G$ is an algebraic group  acting on $X$, we let $\cQ_G(X)^\naive$
denote the dg category of complexes of $G$-equivariant quasicoherent sheaves.
We let $\cQ_G(X)$ denote the localization of this category with respect to acyclic complexes.
We use $\Perf(X) \subset \cQ(X)$ and $\Perf_G(X)\subset \cQ_G(X)$ to denote the full dg subcategories
consisting of \emph{perfect} objects---that is, objects which are quasi-isomorphic to bounded complexes of vector bundles.
If $u:X \to Y$ is a morphism of schemes, we have natural dg functors $u_*:\cQ(X) \to \cQ(Y)$ and
$u^*: \cQ(Y) \to \cQ(X)$.
Note that the functor $u^*$ carries $\Perf(Y)$ to $\Perf(X)$.
Suppose $G$ and $H$ are algebraic groups, $X$ is a scheme with a $G$-action, and $Y$ is a scheme with an $H$-action.
If a morphism $u:X \to Y$ is equivariant with respect to a homomorphism of groups $\phi:G \to H$, then we will often
abuse notation and write $u_*$ and $u^*$ for the equivariant pushforward and pullback functors
$u_*:\cQ_G(X) \to \cQ_H(Y)$ and $u^*:\cQ_H(Y) \to \cQ_G(X)$.

\subsection{Constructible and microlocal geometry}

We refer to \cite{KS} for the microlocal theory of sheaves.
If $X$ is a topological space, let $\Sh(X)$ denote the dg category of bounded chain complexes of sheaves of $\bC$-vector spaces on $X$, localized with respect to acyclic complexes.  If $X$ is a real-analytic manifold, $\Sh_c(X)$ denotes the full subcategory of $\Sh(X)$ of objects whose cohomology sheaves are constructible with respect to a real-analytic
Whitney stratification of $X$.
If $X$ is a (possibly non-compact) real-analytic manifold, then $\Sh_{cc}(X) \subset \Sh_c(X)$ is the full subcategory of objects which have compact support.
We continue to use the phrase ``sheaf'' for an object of $\Sh_{cc}(X)$.

The \emph{standard constructible sheaf} on a submanifold $i: Y\hookrightarrow X$ is defined as the push-forward of the constant sheaf on $Y$, i.e. $i_{*} \bC_Y$ as an object in $Sh_c(X)$. The Verdier duality functor $\cD: Sh_c^\circ(X)\to \Sh_c(X)$ takes $i_{*}\bC_Y$ to the \emph{costandard constructible sheaf} on $X$. We know $\cD(i_{*}\bC_Y)=i_{!}\cD(\bC_Y)=i_{!} \omega_Y$, where $\omega_Y=\cD(\bC_Y)=\bC_Y[\dim Y]$.

We denote the singular support of a complex of sheaves $F$ by $\SS(F) \subset T^*X$.  If $X$ is a real-analytic manifold and $\Lambda \subset T^*X$ is an $\bR_{>0}$-invariant Lagrangian subvariety, then $\Sh_c(X;\Lambda)$ (resp. $\Sh_{cc}(X;\Lambda)$) denotes the full subcategory of $\Sh_c(X)$ (resp. $\Sh_{cc}(X)$)
whose objects have singular support in $\Lambda$.

\subsection{Toric geometry}
\label{sec:toric-geometry}

Let $N \cong \bZ^n$ be a free abelian group, and let $\Sigma$ be a \emph{fan} in $N$ (or in $N_\bR = N \otimes \bR$)
of strongly convex rational polyhedral cones.
We do not necessarily assume that $\Sigma$ satisfies further conditions---e.g. that it is complete, or simplicial.

\subsubsection{Notation}
Given $N$ and $\Sigma$, we fix the following standard notation:
\begin{itemize}
\item  $M := \Hom(N,\bZ)=:N^\vee$ is the dual lattice to $N$.
\item  $N_\bR$ and $M_\bR$ are the real vector spaces spanned by $N$ and $M$, i.e. $N_\bR = N \otimes_\bZ \bR$ and
$M_\bR = M \otimes_\bZ \bR$.
\item  $X_\Sigma$ is the complex toric variety associated to $\Sigma$.  It is naturally equipped with an action of the algebraic torus $T = N \otimes \bC^*$.
\end{itemize}

We also use:
\begin{itemize}
\item  $T_\bR$ denotes the maximal compact subgroup of $T$.  So $T_\bR \cong N_\bR/N\cong \mathrm{U}(1)^n$.
\item  Dually, $T^\vee:=M\otimes \bC^*$ and $T^\vee_\bR\cong M_\bR/M$ is its maximal compact.
\item $\Si(d)$ is the set of $d$-dimensional cones in $\Si$. In particular,
$\Si(1)=\{ \rho_1,\ldots,\rho_r\}$ is the set of rays.
Let $v_i\in N$ be the generator of $\rho_i$, i.e. $\rho_i\cap N = \bZ_{\geq 0} v_i$.
\item Let $\langle \ , \rangle: M_\bR\times N_\bR \to \bR$ denote
the natural pairing.
\item Given a cone $\sigma\in \Sigma$, let
$$
\sigma^\vee=\{ x\in M_\bR\mid \langle x, y \rangle\geq 0 \textup{ for all } y\in N_\bR\}
$$
be the dual cone, and define
$$
\sigma^\perp=\{x\in M_\bR\mid \langle x,y\rangle =0 \textup{ for all } y\in N_\bR\}.
$$
If $\sigma$ is a $d$-dimensional cone then
$\sigma^\perp\subset M_\bR$ is a codimension-$d$ $\bR$-linear subspace.
\end{itemize}

\subsubsection{Equivariant line bundles}
 Let $D_i$ be the $(n-1)$-dimensional $T$ orbit closure
associated to $\rho_i$, so that $D_i$ is a $T$-divisor of $X$.
Any $T$-divisor $D$ of $X$ is of the form $D_\vc= \sum_{i=1}^r c_i D_i$,
where $\vc=(c_1,\ldots,c_r)\in \bZ^r$,
and any $T$-equivariant line bundle on $X$ is of the form $\cL_\vc = \cO_{X_\Si}(D_\vc)$.
If $\cL_\vc$ is ample then
\begin{equation}\label{eqn:polytope-c}
\tri_\vc :=\{  m\in M_\bR\mid \langle m, v_i\rangle \geq - c_i,  i=1,\ldots,r \}
\end{equation}
is a convex polytope in $M_\bR$.

\subsubsection{Orbits}

The $T$-orbits of $X_\Sigma$ can be described using the
structure of the fan.   Given a $d$-dimensional cone $\tau \in\Si$, let
$\tau^\perp$ be the $(n-d)$-dimensional subspace of
$M_\bR$ defined by
$$
\tau^\perp=\{ m\in M_\bR\mid \langle m, y\rangle =0\ \forall y\in \si\}.
$$
Let $N_\tau$ be the rank $d$ sublattice of $N$ generated by $\tau\cap N$,
and let
$$
N(\tau)=N/N_\tau,\quad M(\tau)= \tau^\perp\cap M.
$$
Then $N(\tau)$ and $M(\tau)$ are dual lattices of rank $(n-d)$, and
$$
O_\tau=\Hom(M(\tau),\bC^*)=\Spec\, \bC[M(\tau)]=N(\tau)\otimes \bC^* \cong (\bC^*)^{n-d}
$$
is a $T$-orbit in $X_\Si$. The stabilizer of any point in $O_\tau$ is
$T_\tau := N_\tau\otimes \bC^* \cong (\bC^*)^d$.
In particular, $O_{\{ 0\}}=\Hom(M,\bC^*)= N\otimes \bC^* = T\cong(\bC^*)^n$.

We have a disjoint union of $T$-orbits:
\begin{equation}\label{eqn:Corbit}
X_\Si=\bigcup_{\tau\in \Si} O_\tau,
\end{equation}
which is a $T$-equivariant stratification of $X_\Si$.
Let
$$
X_\tau = \Spec\,\bC[\tau^\vee \cap M] \cong (\bC^*)^{n-d}\times \bC^d
$$
be the affine toric subvariety of $X_\Si$ associated to $\tau$. There
is an inclusion $O_\tau \subset X_\tau$ and a deformation
retraction $r_\tau: X_\tau\to O_\tau$. More explicitly, there
exists a basis $w_1,\ldots, w_d$ of $N_\tau$ such that
$$
\tau = \{ r_1 w_1+ \ldots + r_d w_d \mid r_i\geq 0\}.
$$
Define
$$
\tau^\circ =\{ r_1 w_1+\ldots + r_d w_d \mid r_i >0\}.
$$
Suppose that $y\in N_\tau\cap \tau^\circ$, so that
$y= \sum_{j=1}^d n_j w_j $ where $n_j \in \bZ_{>0}$.
Then the retraction $r_\tau$ is given by
$$
r_\tau(p)=\lim_{t\to -\infty} e^{ty}\cdot p.
$$
Since $T=O_{\{ 0\}}$ is contained in $X_\tau$ for all $\tau\in \Si$, we have
a surjective map $r_\tau: O_{\{ 0\}} \to O_\tau$ which can be identified
with the natural projection $T\to T/T_\tau$.

There is an inclusion $j:\bR^+=\{ e^y\mid y\in \bR \}
\hookrightarrow \bC^*= \{ e^y\mid y\in \bC\} =\bC^*$
and a retraction $r: \bC^* \to \bR^+$ given by $ z\mapsto |z|$. This
induces inclusions
$$
j: O_\tau^+ \stackrel{\mathrm{def}}{=}
N(\tau)\otimes \bR^+ \hookrightarrow O_\tau= N(\tau)\otimes \bC^*
$$
and retractions
$$
r: O_\tau=
N(\tau)\otimes\bC^* \to O_\tau^+=N(\tau)\otimes \bR^+.
$$
In particular, $O_{\{0\}}^+ \cong (\bR^+)^n$ is the image of
the inclusion $\exp: N_\bR \to T$ given by $y\mapsto \exp(y)$.
For each $\tau\in \Si$, we have a surjective map
$$
r_\tau^+:O_{ \{ 0\} }^+ \cong (\bR^+)^n \to  O_\tau^+
\cong (\bR^+)^{n-d}.
$$

Let
\begin{equation}\label{eqn:Rorbit}
(X_\Si)_{\geq 0} = \bigcup_{\tau \in \Si} O_\tau^+.
\end{equation}
Then we have an inclusion
$j: (X_\Si)_{\geq 0}\hookrightarrow X_\Si$ and
a retraction $r: X_\Si\to (X_\Si)_{\geq 0}$.
The retraction $r$
descends to a homeomorphism
$X_\Si/T_\bR \cong (X_\Si)_{\geq 0}$.

\section{Homological Mirror Symmetry for Toric Varieties}
\label{sec:HMS}

In this section, we derive theorems relating the category of coherent
(equivariant) sheaves on $X_\Si$ to the Fukaya category on $T^*T^\vee_\bR$
($T^*M_\bR$). 

\subsection{The Coherent-Constructible Correspondence}

In this subsection, we briefly recall the results of \cite{FLTZ, T}.
The results in \cite{FLTZ,T} hold for toric varieties over an arbitrary commutative,
Noetherian base ring $\mathsf{R}$. Here we state the results
for the case $\mathsf{R}=\bC$. We use the notation in Section \ref{sec:notation}.

Let $X_\Si$ be a toric variety defined by a {\em complete} fan
$\Si\subset N_\bR$. We define
\begin{equation}\label{eqn:LS}
\LS =\bigcup_{\tau\in \Si}(\tau^\perp + M) \times -\tau \subset M_\bR \times N_\bR = T^*M_\bR.
\end{equation}
where $\tau^\perp + M=\{ x+\chi\mid x\in \tau^\perp, \chi\in M\}$.
Then $\LS$ is a Lagrangian subvariety of $T^*M_\bR$.
Let $\bar{\Lambda}_\Si \subset T^*T_\bR^\vee$ be the image of
$\LS$ under the universal covering map
$T^*M_\bR = M_\bR\times N_\bR \to T^*T_\bR^\vee = T_\bR^\vee \times N_\bR$.

\begin{theorem}[equivariant coherent-constructible correspondence {\cite{FLTZ}}] \label{ccc-T}
Let $X_\Si$ be a complete toric variety defined by a finite complete fan $\Si\subset N_\bR$.
Then there is a quasi-equivalence of monoidal dg categories
\begin{equation}\label{eqn:ccc-T}
\kappa: \Perf_T(X_\Sigma)\longrightarrow Sh_{cc}(M_\bR;\LS).
\end{equation}
The functor $\kappa$ sends an equivariant ample line bundle
$\cL_\vc$ on $X_\Si$  to the costandard constructible sheaf $i_!\omega_{\tri_\vc^\circ}$
on $M_\bR$, where $\tri_\vc^\circ$ is the interior of the convex polytope
$\tri_\vc$.
\end{theorem}

\begin{theorem}[nonequivariant coherent-constructible correspondence {\cite{T}}]\label{ccc}
There is a quasi-embedding of monoidal dg categories:
\begin{equation}\label{eqn:ccc}
\bar{\kappa}: \Perf(X_\Sigma) \longrightarrow Sh_c(T_\bR^\vee;\bar{\Lambda}_\Si)
\end{equation}
which makes the following square commute up to natural isomorphism:
\begin{equation}
\begin{CD}
\Perf_T(X_\Si) @>{\kappa}>> Sh_{cc}(M_\bR;\LS)\\
@V{f}VV  @V{p_!}VV\\
\Perf(X_\Si) @>{\bar{\kappa}}>> Sh_c(T_\bR^\vee;\bar{\Lambda}_\Si)
\end{CD}
\end{equation}
where $f$ forgets the equivariant structure, and
$p:M_\bR \to T_\bR^\vee = M_\bR/M$ is the natural projection.
\end{theorem}

\begin{remarks}\label{remark-ccc}
\begin{enumerate}
\item  The monoidal structures in Theorem \ref{ccc-T} and Theorem \ref{ccc} will be discussed
in Section \ref{sec:functoriality} below.
\item  In \cite{FLTZ}, $\kappa$ is defined in terms of certain equivariant quasicoherent
sheaves that arise naturally in the \v{C}ech resolution.
\item In \cite{T}, the third author proved that \eqref{eqn:ccc} is a quasi-equivalence
when $X_\Si$ is a projective, unimodular, zonotopal toric variety.  
We conjecture that \eqref{eqn:ccc} is a quasi-equivalence for any complete toric variety.
\end{enumerate}

\end{remarks}

\subsection{The unwrapped Fukaya category}
\label{unwrapped}
The Fukaya category of the cotangent $T^*X$
of a compact real analytic manifold $X$ was defined
in \cite{NZ} and equated with constructible sheaves on $X$ in
\cite{NZ,N}.  Here we review aspects most
relevant to the present case, including the role of infinity and of standard branes.

Let $X$ be a real analytic manifold equipped with a
Riemannian metric $g$. Let $\pi: T^*X \to X$ be the cotangent bundle
of $X$. Define the closed unit disc bundle to be
$$
D^*X=\{(x,\xi)\in T^*X \mid \|\xi\|\le 1\} \subset T^*X,
$$
and define the unit sphere bundle to be
$$
S^*X=\{(x,\xi)\in T^*X\mid \|\xi\|=1\}=\partial(D^*X)\subset T^*X.
$$
We may think of $D^*X$ as a compactification $\overline{T}^*X$ of $T^*X$ by the following compactification map
\begin{equation}\label{eqn:iota}
\begin{array}{ccccc}
\iota: &T^*X & \to &  D^*X\\
       \  &(x, \xi) &\mapsto & \Bigl(x,\displaystyle{\frac{\xi}{\sqrt{1+\|\xi\|^2}}} \Bigr).
\end{array}
\end{equation}
and we can think of $S^*X$ as $T^\infty X$ because it is the ``infinity'' part of $T^*X$ under this compactification.
If $L$ is a Lagrangian submanifold of $T^*X$ we write $L^\infty$ for $\overline{\iota(L)}\cap S^*X$, the
part of $L$ at infinity in the fibers.

In the present case, $M_\bR$ is noncompact.  This is only a minor complication,
as we will require all
Lagrangian branes $L$ to have compact {\em horizontal support}, i.e.,  $\overline{\pi(L)}$
is compact.  Define the flat metric $g$ on $T^*M_\bR = M_\bR \times N_\bR$ by declaring a $\bZ$ basis
$\{ e_1,\ldots, e_n\}$ of $N\subset N_\bR$ to be orthonormal,
and likewise for the dual basis $\{ e_1^*,\ldots, e_n^*\}$ of $M.$  Then $\overline{\pi(L)}$ is bounded.
%
%
We require as well that
the usual other conditions of Lagrangian branes are satisfied:  that is,
$L$ must be an exact Lagrangian submanifold of $T^*M_\bR$;
$\overline{\iota(L)}$ is a $\cC$-set of $D^*X$;\footnote{See Section
\ref{sec:review-Cset} for a brief review of analytic-geometric categories, including
definitions of $\cC$-sets and $\cC$-maps.}
and $L$ is equipped with the data of a vector bundle with
flat connection, a brane structure and a tame perturbation
(see \cite{NZ}). Under these conditions, morphisms are well-defined for the
following reason.  If $L=(L_1,...,L_k)$ is a finite collection of Lagrangian objects
with compact horizontal support, then there exists
a sublattice $\Xi\subset M$ of finite index $d$ such that the union of the supports
of the $L_i$ are contained in a single fundamental domain:  then all morphisms can be
computed in the cotangent of
the compact torus $M_\bR/\Xi$ (a degree $d$ cover of
the torus $T_\bR^\vee=M_\bR/M$) and lifted to $T^*M_\bR =N_\bR\times M_\bR$---see
\cite[Section 5.3]{N} for details.\footnote{The condition
of compact horizontal support can be dropped for a single given
object, as one can define the Yoneda image by analyzing hom's
against objects with compact horizontal support---see \cite{N2}.}
Holomorphicity is preserved by the lift
since the quotient by $\Xi$ is a local isomorphism of the
K\"{a}hler structure. The triangulated envelope of the Fukaya
$A_\infty$-category of all such branes is denoted by $Fuk(T^*M_\bR)$.\footnote{The triangulated envelope of any $A_\infty$-category is unique up to an exact quasi-equivalence.}

Let $\Lambda \subset T^*X$ be a conical Lagrangian subset.  The $A_\infty$-category
generated by Lagrangian branes $L$ with $L^\infty \subset \Lambda^\infty$ is
denoted by $Fuk(T^*X;\Lambda).$  Here we will mainly be concerned with
$Fuk(T^*M_\bR;\LS)$, where $\LS$ is given in \eqref{eqn:LS}.

\subsection{Microlocalization}

Recall that if $i:Y\hookrightarrow X$ is the inclusion of an analytic submanifold
in a compact, real analytic manifold $X$ then $i_* \bC_Y$ is the standard object in $Sh_{c}(X)$
associated to $Y$, and under microlocalization, the {\em standard brane} $\mu(i_*\bC_Y)$ is defined
by the {\em standard Lagrangian}
$L_{Y,*}\subset T^* X$ given by the fiberwise sum
$$
L_{Y,*}= T^*_Y X +\Gamma_{d f},
$$
where $f = \log m$ and $m$ is a nonnegative $\cC$-function
$m: X \to \bR$ that vanishes precisely on the boundary $\pa Y\subset X$.
Here $T^*_Y X$ is the conormal bundle of $Y$ in $X$,
and $\Gamma_{df} \subset T^*Y \cong T^*X/T^*_Y X$
is the graph of $df$.
There is a canonical brane structure on this Lagrangian (Section 5.3 of \cite{NZ}).
We let $L_{Y,m,*}$ denote the standard Lagrangian defined by a particular
choice of $m$. Two different choices $m_1$, $m_2$ give rise
to isomorphic objects: $L_{Y,m_1,*}\cong L_{Y,m_2,*}$ as objects
in $Fuk(T^*X)$.


Let $\alpha$ be a diffeomorphism on $M_\bR\times N_\bR$ given by $\alpha(x,y)=(x,-y)$.
A \emph{costandard brane} (\emph{costandard Lagrangian}) $L$ is a brane (Lagrangian)
such that $\alpha(L)$ is a standard brane (Lagrangian).
Microlocalization $\mu$
also takes the costandard constructible sheaf $i_! \omega_Y$ to the costandard brane $L_{Y,!}:=
T^*_Y X - \Gamma_{df}.$
We summarize these results as a theorem.
\begin{theorem} [\cite{NZ,N}]\label{nadler-zaslow}
There is a quasi-equivalence of $A_\infty$-categories
$$
\mu: Sh_{cc}(M_\bR;\LS)\to Fuk(T^*M_\bR;\LS).
$$
For any analytic submanifold $Y\subset M_\bR$,
$\mu$  takes the standard constructible sheaf $i_*\bC_Y$ to the standard brane
$L_{Y,*}$ of $Y$, and takes the costandard constructible sheaf $i_!(\omega_Y)$
to the costandard brane $L_{Y,!}.$
The functor $\mu$ admits a quasi-inverse $\mu^{-1}: Fuk(T^*M_\bR;\LS)\to Sh_{cc}(M_\bR;\LS)$.

Similarly, we have a quasi-equivalence of $A_\infty$-categories
\begin{equation}\label{eqn:mu}
\bar{\mu}: Sh_c(T^\vee_\bR;\bar{\Lambda}_{\Si})
\stackrel{\cong}{\longrightarrow}
Fuk(T^* T^\vee_\bR;\bar{\Lambda}_\Si).
\end{equation}
\end{theorem}

\subsection{Functoriality and monoidal structure}
\label{sec:functoriality}

The functoriality of the functor $\kappa$ is proven in \cite{FLTZ}.
We extend this functoriality involving the Fukaya category.
This is simply a combination of the results of \cite{FLTZ} and of \cite[Section 5]{N}.

We first review some general results in \cite[Section 5]{N}.
Given two real analytic manifolds $X_0$, $X_1$, let
$p_0:X_0\times X_1\to X_0$ and $p_1:X_0 \times X_1\to X_1$
be projections. For a real analytic manifold $Y$,
let $\mu_Y: Sh_c(Y)\to Fuk(T^*Y)$
be the microlocalization functor, and let
$\alpha_Y:Fuk(T^*Y)\to Fuk(T^*Y)^\circ$
be the brane duality functor (induced by multiplication
by $-1$ on cotangent vectors). Let
$$
\cY_\ell :Fuk(T^*X_1)\to mod_l(T^*X_1)^\circ,\quad
P\mapsto hom_{Fuk(T^*X_1)}(P,-)
$$
be the Yoneda embedding for left
$A^\infty$-modules over $Fuk(T^*X_1)$.

 An object $\cK$ of $\Sh_c(X_0\times X_1)$
defines a functor
\begin{equation}
\Phi_{\cK !}: \Sh_c(X_0)\to \Sh_c(X_1),\quad
\cF \mapsto p_{1!}(\cK\otimes p_0^* \cF).
\end{equation}
An object $L$ of $Fuk(T^* X_0\times T^*X_1)$ defines a functor
\begin{equation}
\begin{aligned}
\tilde{\Psi}_{L!}: & Fuk(T^* X_0) \to mod_\ell (Fuk(T^*X_1))^\circ,\\
& P\mapsto  hom_{Fuk(T^*X_0\times T^*X_1)}(L, \alpha_{X_0}(P)\times -)
\end{aligned}
\end{equation}
The following is a special case of  \cite[Proposition 5.3.1]{N}.
\begin{theorem}\label{thm:integral-transform}
Consider an object $\cK$ of $\Sh_c(X_0\times X_1)$, and its
microlocalization $L=\mu_{X_0\times X_1}(\cK)$. Then
there is a quasi-isomorphism
$$
\cY_\ell\circ \mu_{X_1}\circ \Phi_{\cK!}\simeq
\tilde{\Psi}_{L!} \circ \mu_{X_0}.
$$
Therefore the functor $\tilde{\Psi}_{L!}$ is represented by
$$
\Psi_{L!}:= \mu_{X_1}\circ \Phi_{\cK!}\circ \mu_{X_0}^{-1} :
Fuk(T^*X_0)\to Fuk(T^*X_1).
$$
\end{theorem}

\begin{example}\label{Fuk-functorial}
Let $v:X_0\to X_1$ be a smooth map, and let
$$
\Gamma_v=\{(x_0,x_1)\in X_0\times X_1\mid x_1=v(x_0)\}
$$
be the graph of $v$. Let $\cK=\bC_{\Gamma_v}$
be the constant sheaf on $\Gamma_v$, and
let $L_v= \mu_{X_1\times X_2}(\cK)$. Then
$$
\Phi_{\cK!}=v_!,\quad L_v \simeq T_{\Gamma_v}^*(X_0\times X_1)
$$
where $T_{\Gamma_v}^*(X_0\times X_1)$ is the conormal bundle
of $\Gamma_v$ in $X_0\times X_1$. Define $\Psi_{L_v!}$
as in Theorem \ref{thm:integral-transform}. Then
there is a quasi-isomorphism
\begin{equation}\label{eqn:Fuk-functorial}
\Psi_{L_v !} \circ \mu_{X_0}\simeq \mu_{X_1}\circ v_!.
\end{equation}
\end{example}

For two toric varieties $X_1=X_{\Si_1}$ and $X_2=X_{\Si_2}$ and a fan-preserving map $f: N_1\to N_2$,
let $v:M_{2,\bR}\to M_{1,\bR}$  and $u:X_1 \to X_2$
be the induced map of vector spaces and varieties (see \cite{FLTZ}).
As a special case of Example \ref{Fuk-functorial},
define  $L_v := T^*_{\Gamma_v}(M_{2,\bR}\times M_{1,\bR})$,
which is a Lagrangian subspace
of $T^*M_{2,\bR}\times T^* M_{1,\bR}=
M_{\bR,2}\times N_{\bR,2}\times M_{\bR,1}\times N_{\bR,1} $.
Combining \eqref{eqn:Fuk-functorial} with
Theorem \ref{thm:integral-transform} and results in 
\cite[Section 3]{FLTZ}, we come to a larger diagram:
\begin{theorem}\label{thm:Fuk-functoriality}
 For two complete toric varieties $X_1=X_{\Si_1}$ and
$X_2=X_{\Si_2}$  and a fan-preserving map $f: N_1\to N_2$, where $f$ is
injective, and associated maps $f\otimes 1_{\bC^*}:T_1\to T_2$,
$u: X_1\to X_2$, $v: M_{2,\bR}\to M_{1,\bR}$, the following diagram
commutes up to a quasi-isomorphism.
$$\xymatrix{
 \Perf_{T_2}(X_2) \ar[d]_{u^*}  \ar[r]^{\kappa_2}  &Sh_{cc}(M_{2,\bR}; \Lambda_{\Si_2}) \ar[r]^-{\mu_{M_{2,\bR}}}  \ar[d]^{v_!}       & Fuk(T^*M_{2,\bR};\Lambda_{\Si_2})\ar[d]^{\Psi_{L_v!}}&\\
\Perf_{T_1}(X_1) \ar[r]_{\kappa_1}                           & Sh_{cc}(M_{1,\bR}; \Lambda_{\Si_1})
\ar[r]_-{\mu_{M_{1,\bR}}}                    & Fuk(T^*M_{1,\bR};\Lambda_{\Si_1}).
}$$
\end{theorem}

\begin{example}[a product structure on the Fukaya category]\label{Fuk-tensorial}
This example is a special case of Example \ref{Fuk-functorial}.

Let $G$ be a Lie group, and let $v: G \times G\to G$ be
the multiplication: $v(g_1, g_2) = g_1\cdot g_2$.
Then $L_v$ is an object in $Fuk(T^*(G\times G)\times T^*G)$ and defines
a functor $\Psi_{L_v !}: Fuk(T^*G\times T^*G)\to Fuk(T^*G)$.
We define the product $L_1 \diamond L_2 $ of two objects $L_1$ and $L_2$
of $Fuk(T^* G)$ by the formula
\begin{equation}\label{eqn:diamond}
L_1 \diamond L_2 := \Psi_{L_v!}(L_1\times L_2).
\end{equation}
\end{example}

\begin{proposition}[the microlocalization intertwines the product structures]
\label{tensor}
Let $G$ be a Lie group.
The microlocalization functor $\mu_G: Sh_{cc}(G)\stackrel{\sim}{\to} Fuk(T^*G)$ intertwines the monoidal product on  $Sh_{cc}(G)$ given by the convolution, and the product structure on $Fuk(T^*G)$ given by the product $\diamond$ defined by \eqref{eqn:diamond}, up to a quasi-isomorphism: i.e. the functors $\mu_{G}(-\star-)$ and $\mu_G(-)\diamond \mu_G(-)$ are quasi-isomorphic in the category of $A_\infty$-functors from $Sh_{cc}(G)\times Sh_{cc}(G)$ to $Fuk(T^*G)$.
\end{proposition}
\begin{proof} Recall that convolution product
$F_1\star F_2$ of two objects  $F_1$ and $F_2$ of $\Sh_{cc}(G)$ is defined by
$F_1\star F_2 = v_!(F_1\boxtimes F_2)$. So
\begin{eqnarray*}
\mu_G (F_1\star F_2) &=& \mu_G\circ v_!(F_1 \boxtimes F_2) \cong  \Psi_{L_v!} \circ \mu_{G\times G}(F_1\boxtimes F_2)\\
&=&  \Psi_{L_v!}(\mu_G(F_1)\times \mu_G(F_2)) = \mu_G(F_1)\diamond \mu_G(F_2)
\end{eqnarray*}
\end{proof}

\subsection{Equivariant and nonequivariant HMS for toric varieties}
Let $\tau=\mu\circ \kappa$ and let $\bar{\tau}= \bar{\mu}\circ
\bar{\kappa}$. Notice that the convolution product of costandard
sheaves $i_{1!} \omega_{\Delta_{\vc_1}}$ and ${i_{2!}}
\omega_{\Delta_{\vc_2}}$ is $i_! \omega_{\Delta_{\vc_1+\vc_2}}$,
where $c_1$ and $c_2$ determine two equivariant ample line bundles
on $X_\Si$, and $i_1$, $i_2$ and $i$ are corresponding embeddings of
polytopes. Since costandard sheaves over convex polytopes of ample
line bundles generate the category $Sh_{cc}(M_\bR;\LS)$, as shown in
\cite{FLTZ}, the subcategory
$Sh_{cc}(M_\bR;\LS)$ of $Sh_{cc}(M_\bR;\LS)$ is closed under the convolution product. 
By results in \cite{T}, the subcategory $Sh_c(T_\bR^\vee,\bar{\Lambda}_\Si)$
of $Sh_c(T_\bR^\vee)$ is closed under the convolution product.
Combining Theorem \ref{ccc-T}, Theorem \ref{ccc}, Theorem
\ref{nadler-zaslow}, and Proposition \ref{tensor}, we obtain:
\begin{theorem}\label{pre-hms}
Let $X_\Si$ be a complete toric variety defined by a finite complete fan
$\Si\subset N_\bR$. Then there is an quasi-equivalence of $A_\infty$
categories
\begin{equation}\label{eqn:tau}
\tau: \Perf_T(X_\Si)\stackrel{\cong}{\longrightarrow} Fuk(T^*M_\bR;\LS).
\end{equation}
There is an quasi-embedding of $A_\infty$ categories
\begin{equation}\label{eqn:bar-tau}
\bar{\tau}: \Perf(X_\Si) \to  Fuk(T^*T_\bR^\vee;\bar{\Lambda}_\Si).
\end{equation}
The functors $\tau$ and $\bar \tau$ intertwine the corresponding product structures in the Fukaya categories and the monoidal products in the dg category of perfect sheaves. More precisely, there are quasi-isomorphisms
$$\tau(-\otimes-)\cong \tau(-)\diamond \tau(-),\ \bar\tau(-\otimes -) \cong \bar\tau(-) \diamond \bar\tau(-).$$
\end{theorem}

Let $DCoh_T(X_\Si)$ be the bounded derived category of $T$-equivariant coherent
sheaves on $X_\Si$, and let $DCoh(X_\Si)$ be the bounded derived category
of coherent sheaves on $X_\Si$. When $X_\Si$ is nonsingular, we have
$D\Perf_T(X_\Si)= DCoh_T(X_\Si)$ and $D\Perf(X_\Si)= DCoh(X_\Si)$.
Taking $H^0$ of \eqref{eqn:tau} and \eqref{eqn:bar-tau}, we obtain
the following Corollary \ref{H-tau} and Corollary \ref{H-bar-tau}, respectively.

\begin{corollary}[Equivariant homological mirror symmetry
of toric varieties]\label{H-tau}
Let $X_\Si$ be a nonsingular complete toric variety defined by a finite nonsingular complete fan
$\Si\subset N_\bR$. Then there is an equivalence of tensor triangulated
categories
\begin{equation}\label{eqn:H-tau}
H(\tau): D Coh_T(X_\Si) \stackrel{\cong}{\longrightarrow} DFuk(T^*M_\bR;\LS)
\end{equation}
\end{corollary}

\begin{corollary}[Homological mirror symmetry of toric varieties] \label{H-bar-tau}
Let $X_\Si$ be a nonsingular complete toric variety defined by a finite nonsingular complete fan
$\Si\subset N_\bR$.
There is an embedding of tensor triangulated categories
\begin{equation}\label{eqn:H-bar-tau}
H(\bar{\tau}): D Coh(X_\Si)\longrightarrow DFuk(T^*T_\bR^\vee;\bar{\Lambda}_\Si).
\end{equation}
\end{corollary}
We conjecture that \eqref{eqn:H-bar-tau} is an equivalence.  This is
proven for $\Si$ a complete, unimodular hyperplane arrangement in \cite{T}.

\section{T-duality}
\label{sec:tduality}

In this section, we perform an equivariant version of T-duality.
Let $X_\Si$ be an $n$-dimensional nonsingular projective toric variety (so that it is a compact toric manifold).
Then $T\cong (\bC^*)^n$ and its maximal compact subgroup
$T_\bR\cong U(1)^n$ acts on $X_\Si$. From a $T$-equivariant line
bundle $\cL_\vc$ on $X$ together with a $T_\bR$-invariant
hermitian metric $h$, we construct a Lagrangian submanifold
$\bL_{\vc,h}$ of $T^*M_\bR$.
We relate $\bL_{\vc,h}$ to classical objects in symplectic  geometry.

%

\subsection{Construction of the T-dual Lagrangian}
\label{sec:construct}
Let $X= X_\Si$ be a smooth projective toric variety defined
by a fan $\Si\subset N_\bR$, and let
$\rho_1,\ldots,\rho_r$ be the 1-dimensional cones
in $\Si$ and $D_1,\ldots,D_r$ the associated $T$-divisors,
as in Section \ref{sec:toric-geometry}.

There exists $s_i\in H^0(X, \cO_X(D_i))$, unique
up to multiplication by a constant scalar in $\bC^*$, such that the zero locus of
$s_i$ is exactly $D_i$.
$$
X_{\{0\}} = X\setminus \cup_{i=1}^r  D_i =\Spec \bC[M]\cong (\bC^*)^n.
$$
is the unique open orbit of the $T$-action.

The meromorphic section $s_\vc:=\prod_{i=1}^r s_i^{c_i}$
of $\cL_\vc =\cO_X(D_\vc)$ is defined up to multiplication by a constant scalar in $\bC^*$.
The restriction of $s_\vc$ to $X_{\{0\}}$ is a holomorphic frame of $\cL_\vc$ on
the Zariski open subset $X_{ \{ 0\} } \subset X$.

We now choose a $T_\bR$-invariant, real analytic hermitian metric
$h$ on $\cL_\vc$.
Let $\nabla_{\vc, h}$ be the unique connection on $\cL_\vc$
determined by the  holomorphic structure on $\cL_\vc$ and the hermitian
metric $h$. The connection 1-form of $\nabla_h$ with respect to
the unitary frame $s_\vc/\|s_\vc\|_h$ of $\cL_\vc|_{X_{\{ 0\}}}$ is the
following purely imaginary, real analytic 1-form.
$$
\alpha = -2\sqrt{-1}\mathrm{Im}(\dbar\log \| s_\vc \|_h).
$$

Note that $\alpha$ is invariant if we replace $s_\vc$ and $h$
by $\lambda s_\vc$ and $\rho h$ respectively, where
$\lambda\in \bC^*$ and $\rho\in (0,\infty)$ are constants.

We now introduce coordinates on $X_{\{0\}}\cong T$
(the identification depends on the choice of a point
in $X_{\{0\}}$).
The universal cover of $T$ can be canonically identified with
$N\otimes \bC=N_\bR\times N_\bR$. Let $\{ e_1,\ldots,e_n\}$ be a
$\bZ$-basis of the lattice $N$, and let
$\{e_1^*,\ldots,e_n^*\}$ be a dual $\bZ$-basis of the dual lattice $M$.
A vector in $N\otimes \bC$ is of the form $\sum_{j=1}^n\frac{y_j +\sqrt{-1}\theta_j}{2\pi} e_j$
where $y_j,\theta_j\in \bR$.  A vector in $M_\bR$ is of the
form $\sum_{j=1}^n \frac{\gamma_j}{2\pi} e_j^*$, where $\gamma_j\in \bR$. Then
$y_j+ \sqrt{-1}\theta_j$ are complex coordinates
on $N_\bC$, and $\gamma_j,$ $y_j$ are Darboux coordinates on
$T^*M_\bR =M_\bR \times N_\bR$. The symplectic form on $M_\bR\times N_\bR $
is
$$
\omega^\vee=\sum_{j=1}^n dy_j\wedge d\gamma_i
$$
which descends to a symplectic form on
$(M_\bR/M)\times N_\bR \cong T^*(T_\bR^\vee)$.
Note that $M\subset M_\bR$ is given by $\gamma_j\in 2\pi\bZ$ and
$N\subset N_\bR$ is given by $\theta_j\in 2\pi\bZ$.
Let $r_j = e^{y_j}$, so that the coordinates on  $T$ are
$e^{y_j+ \sqrt{-1}\theta_j}= r_j e^{\sqrt{-1}\theta_j}$,
$j=1,\ldots,n$.

The function $\| s\|_h$ is $T_\bR$-invariant, so it
depends on $r_i$ ($y_i$) but not on $\theta_i$.
We have
\begin{eqnarray*}
\sqrt{-1}\alpha &=&  2 \mathrm{Im} (\dbar \log \| s_\vc\|_h)
=\mathrm{Im} \Bigl(\sum_{j=1}^n \bigl(\frac{\pa}{\pa r_j} \log \| s_\vc \|_h\bigr)\cdot
(dr_j - \sqrt{-1}r_j d\theta_j) \Bigr) \\
&=& -\sum_{j=1}^n \bigl(\frac{\pa}{\pa y_j} \log \| s_\vc \|_h\bigr)  d\theta_j.
\end{eqnarray*}

Let $y=(y_1,\ldots,y_n)$, and let $f_{\vc, h} (y)  = -\log\| s_\vc \|_h$.
Then $f_{\vc, h}(y)$ is a real analytic function in $y$, and
\begin{equation}
\sqrt{-1}\alpha = \sum_{j=1}^n \frac{\partial f_{\vc,h}}{\partial y_j} (y) d\theta_j.
\end{equation}

We now T-dualize following \cite{AP}.
Specifically, the data of a Lagrangian
section of the dual torus fibration $T_\bR ^\vee\times N_\bR \to N_\bR$
(projection to the second factor)
is equated with a $T_\bR$-invariant $U(1)$-connection
on the torus fibration $p_2:T_\bR\times N_\bR \to N_\bR$ (projection
to the second factor). The restriction
of $\alpha$ to a fiber $p_2^{-1}(y)\cong T_\bR$
is a harmonic 1-form on the torus $p_2^{-1}(y)$, which can be viewed as
an element in $H^1(T_\bR;\bR)\cong M_\bR$,
the universal cover of  the dual torus $T_\bR^\vee=M_\bR/M$ of $T_\bR$.
Let $\bL_{\vc,h}\subset M_\bR\times N_\bR$
be the graph of the map $N_\bR \to M_\bR$ defined by
$y\mapsto \sqrt{-1}\alpha\Bigr|_{p_2^{-1}(y)}$.
In terms of the coordinates $\gamma_j$ on $M_\bR$
and $y_j$ on $N_\bR$,
$\bL_{\vc,h}$ is given by
$$
\frac{\gamma_j}{2\pi} = \frac{\partial f_{\vc,h}}{\partial  y_j} (y), \quad j=1,\ldots, n.
$$
Since $N_\bR$ and $M_\bR$ are dual real vector spaces, we have
$M_\bR\times N_\bR \cong T^*N_\bR \cong T^*M_\bR$. Moreover,
the canonical symplectic forms on $T^* N_\bR$ and
$T^* M_\bR$ are
$$
-\omega^\vee = \sum_{j=1}^n d\gamma_j\wedge dy_j,\quad
\omega^\vee= \sum_{j=1}^n  dy_j \wedge d\gamma_j.
$$
$f_h$  is a  real analytic function on $N_\bR$.
The submanifold $\bL_{\vc, h}$ is the graph of $df_h$ in $T^* N_\bR$, so it is
a real analytic Lagrangian submanifold of
$(T^*N_\bR, -\omega^\vee)$ and
of $(T^*M_\bR, \omega^\vee)$.
Let $\bar{\bL}_{\vc,h} \subset T^* T_\bR^\vee=(M_\bR/M)\times N_\bR$
be the image of $\bL_{\vc,h}$ under the projection $M_\bR\times N_\bR \to (M_\bR/M)\times N_\bR$.
Then $\bar{\bL}_{\vc,h}$ is a real analytic Lagrangian submanifold
in $T^*T_\bR^\vee= T_\bR^\vee\times N_\bR$, and is the graph of a map
$N_\bR\to T_\bR^\vee$. Both $\bL_{\vc,h}$ and $\bar{\bL}_{\vc,h}$
are diffeomorphic to $N_\bR \cong \bR^n$, so they are exact
Lagrangian submanifolds.

Suppose that $D_{\vc'} - D_{\vc}$ is a principal divisor. Then
$\cL_{\vc'}$ and $\cL_{\vc}$ are the same holomorphic  line bundle
equipped with possibly different $T$-equivariant structures, so we may choose
the same hermitian metric $h$ on $\cL_{\vc'}$ and on $\cL_{\vc}$. We have
$s_{\vc'} = s_{\vc} \prod_{j=1}^n t_j^{m_j}$
for some $(m_1,\ldots, m_n)\in \bZ^n$, so
$$
f_{\vc',h} = f_{\vc,h} - \sum_{j=1}^n m_j y_j,\quad
\frac{\partial f_{\vc',h}}{\partial y_j}(y)
= \frac{\partial f_{\vc,h}}{\partial y_j}(y)-  m_j.
$$
Therefore $\bar{\bL}_{\vc',h} =\bar{\bL}_{\vc,h}$.

\subsection{Relations with the equivariant first Chern form and the moment map }
\label{sec:moment}

\subsubsection{Equivariantly closed 2-forms and  moment maps of presymplectic forms}

We recall some definitions from \cite{AB} and \cite{KT}.

The real vector space $N_\bR$ can be identified with the Lie algebra of
the compact torus $T_\bR$, with a basis $\{e_1,\ldots,e_n\}$, and
$M_\bR$ is the dual real vector space, with the dual basis
$\{e_1^*,\ldots,e_n^*\}$. Let $X_j$ be the vector field on $X$
associated to $e_j\in N_\bR$. An equivariant 2-form on $X$ is
of the form
$$
\omega^{\#}= \omega + \sum_{j=1}^n \phi_j e_j^*
$$
where $\omega$ is a $T_\bR$-invariant 2-form on $X$ and $\phi_j$ are $T_\bR$-invariant
functions on $X$. An equivariant 2-form $\omega^{\#}$ is {\em equivariantly closed} if
\begin{equation}\label{eqn:closed}
d\omega =0
\end{equation}
and
\begin{equation}\label{eqn:moment}
i_{X_j}\omega + d\phi_j=0, \quad j=1,\ldots,n.
\end{equation}
In this case, the closed 2-form $\omega$ represents
a cohomology class $[\omega]\in H^2(X;\bR)$, and the equivariantly closed
2-form represents an equivariant cohomology class $[\omega^{\#}]\in H_T^2(X;\bR)$,
and we say $\omega^\#$ (resp. $[\omega^\#]$) is an equivariant lifting
of $\omega$ (resp. $[\omega]$).

In the terminology of \cite{KT}, \eqref{eqn:closed} says that
$\omega$ is a presymplectic form (which is by definition a $T_\bR$-invariant closed 2-form),
and \eqref{eqn:moment} says that $\Phi=\sum_{j=1}^n \phi_j e_j^*:X\to M_\bR$ is a moment map
of the $T_\bR$-action with respect to the presymplectic form $\omega$.
When $\omega$ is nondegenerate, $\omega$ is a symplectic form, and $\Phi$ is a moment
map of the $T_\bR$-action on the symplectic manifold $(X,\omega)$.

\subsubsection{The Equivariant first Chern form}
We now return to the construction in Section \ref{sec:construct}.
Let $F_h$ be the curvature 2-form of the connection $\nabla_h$. Then
$$
c_1(\cL_\vc,\nabla_h)=\frac{\sqrt{-1}}{2\pi} F_h
$$
is a closed, real, $T_\bR$-invariant, real analytic 2-form which
represents the first Chern class
$c_1(\cL_\vc)\in H^2(X;\bR)$. The closed 2-form
$c_1(\cL_\vc,\nabla_h)$ is known as the first Chern {\em form}
defined by the connection $\nabla_h$; it depends
on the underlying holomorphic line bundle and
the hermitian metric $h$, but not on the
equivariant structure.

The section $s_\vc$ determines an equivariant lifting
$c_1(\cL_\vc,\nabla_h,s_\vc)$ of the first Chern form $c_1(\cL_\vc,\nabla_h)$.
More explicitly,
$$
c_1(\cL_\vc,\nabla_h,s_\vc) = \frac{1}{2\pi}(\sqrt{-1}F_h+\sum_{i=1}^n \phi_j e_j^*)
$$
where $\phi_1,\ldots,\phi_n$ are $T_\bR$-invariant, real-valued functions
on $X$. On the open set $X_{\{0\}}\cong (\bC^*)^n$, we have
$X_j=\frac{\partial}{\partial \theta_j}$, and
$$
\sqrt{-1}F_h = \sqrt{-1} d\alpha
=\sum_{j=1}^n d(\frac{\partial f_{\vc,h}}{\partial y_j}) \wedge d\theta_j,
\quad \quad \phi_j =\frac{\partial f_{\vc,h}}{\partial y_j}(y).
$$
The equivariantly closed 2-form $c_1(\cL_\vc,\nabla_h,s_\vc)$ represents
the equivariant first Chern class $(c_1)_T(\cL_\vc)\in H^2_T(X;\bR)$;
we call $c_1(\cL_\vc,\nabla_h,s_\vc)$ the equivariant first
Chern {\em form} defined by $\nabla_h$ and $s_\vc$.

\subsubsection{The Moment Map}
The real analytic map $\Phi_{\vc,h} = \sum_{j=1}^n \phi_j e_j^*:X\to M_\bR$ is a moment map
of the presymplectic form $\omega_h:=\sqrt{-1}F_h$. On $X_{\{0\}}$ it is given by
$$
\Phi_{\vc,h}(y,\theta)= \sum_{j=1}^n \frac{\partial f_{\vc,h}}{\partial y_j}(y)e_j^*.
$$
Define new coordinates
$x_j=\frac{\gamma_j}{2\pi}$ on $M_\bR$, so that $M\subset M_\bR$ is given by
$x_j\in \bZ$. Then the T-dual Lagrangian $\bL_{\vc,h}$ constructed
in Section \ref{sec:construct} can be written as
$$
\bL_{\vc,h}= \{ (x,y)\in M_\bR\times N_\bR\mid x = \Phi_{\vc,h}\circ j_0(y) \}
$$
where $j_0:N_\bR\to X_\Si$ is a composition of inclusions:
$$
N_\bR \stackrel{\exp}{\cong} N\otimes \bR^+\cong (\bR^+)^n
\hookrightarrow (\bC^*)^n \cong N\otimes \bC^* = T =X_\Si -\cup_{i=1}^r D_i
\hookrightarrow X_\Si.
$$
We also have
$$
\Phi_{\vc,h} \circ j_0(N_\bR)= \Phi_h(X_\Si -\cup_{i=1}^r D_i).
$$
The image of $\Phi_{\vc,h}: X\to M_\bR$ is  a {\em twisted polytope}
in the sense of \cite{KT}.

\subsection{T-dual Lagrangians of ample and anti-ample line bundles}
\label{sec:ample-inverse}

When $\cL_\vc$ is ample, we may choose $h$ such that $\omega_h$ is
a symplectic form. Then $\Phi_{\vc,h}: X\to M_\bR$ is
the moment map of the $T_\bR$-action on the symplectic
manifold $(X,\omega_h)$. The
image of the moment map $\Phi_{\vc,h}$ is the convex
polytope $\tri_\vc$ defined by \eqref{eqn:polytope-c}.
Note that the moment map $\Phi_{\vc,h}$ depends on both $\vc$ and $h$, but
the moment polytope $\tri_\vc=\Phi_{\vc,h}(X)$ depends on $\vc$ but not on $h$.
$\Phi_{\vc,h}$ restricts to a homeomorphism
$X_{\geq 0}\to \tri_\vc$, and $\Phi_h\circ j_0: N_\bR \to M_\bR$
maps $N_\bR$ diffeomorphically to $\tri_\vc^\circ$, the
interior of the moment polytope $\tri_\vc \subset M_\bR$.
Let $\Psi_{\vc,h}: N_\bR \to \tri_\vc^\circ$ be this diffeomorphism.
Then $\bL_{\vc,h}$ can be rewritten as a graph over $\tri_\vc^\circ$:
$$
\bL_{\vc,h}=\{ (x,\Psi_{\vc,h}^{-1}(x)) \mid x \in \tri_\vc^\circ \}
\subset \tri_\vc^\circ \times N_\bR = T^* \tri_\vc^\circ \subset T^* M_\bR.
$$

There exists a real analytic function $f^*_{\vc,h}:\tri_\vc^\circ \to \bR$,
unique up to addition of a constant $r\in \bR$, such that
$\Psi_{\vc,h}^{-1}(x)= df^*_{\vc,h}(x)$.  Indeed $f_{\vc,h}^*: \tri_\vc^\circ \to\bR$ can
be chosen to be the {\em Legendre transform} of $f_h:N_\bR\to \bR$.
More explicitly, let $\langle \ ,\ \rangle: M_\bR\times N_\bR\to \bR$ be the natural
pairing. Then
\begin{equation} \label{eqn:legendre}
f^*_{\vc,h}(x) = \sup_{y\in N_\bR}
(\langle x, y\rangle - f_{\vc,h}(y)),
\quad  x\in \tri^\circ_\vc.
\end{equation}

We now consider the equivariant anti-ample line bundle $\cL_\vc^{-1}=\cL_{-\vc}$
equipped with the $T_\bR$-invariant, real analytic hermitian
metric $h^{-1}$. Then $\Phi_{-\vc, h^{-1}}= -\Phi_{\vc,h}$, so
$$
\tri_{-\vc} :=\Phi_{-\vc,h^{-1}}(X)= -\tri_\vc
=\{ m\in M_\bR \mid \langle m, v_i \rangle \leq c_i, i=1,\dots,r\},
$$
and
$$
\bL_{-\vc,h^{-1}}=\{ (-\Psi_{\vc,h}(y),y)\mid y\in N_\bR \}
=\{ (x,\Psi_{\vc,h}^{-1}(-x))\mid x\in \Delta_{-\vc}^\circ \}.
$$
Define a map $\beta: M_\bR\times N_\bR\to M_\bR\times N_\bR$ by $\beta(x,y)=(-x, y)$.
It is easy to see that $\bL_{-\vc, h^{-1}}=\beta(\bL_{\vc, h})$.

\section{T-dual Lagrangians as Objects in the Fukaya Category}
\label{sec:fukaya}

The goal of this section is to prove Theorem \ref{main-III}.
Let $X_\Si$ be a smooth projective toric variety
defined by a fan $\Si\subset N_\bR$.
Let $\cL_\vc$ be an equivariant ample line bundle
on $X_\Si$, and let $\bL_{\vc,h}$ and $\bL_{-\vc,h^{-1}}$
be as in Section \ref{sec:ample-inverse}.
In Section \ref{sec:compactified}, we
prove that $\bL_{\vc,h}$ and $\bL_{-\vc,h^{-1}}$ are objects in $Fuk(T^*M_\bR;\LS)$.
In Section \ref{sec:tlagstd}, we prove that (see Theorem \ref{mainthm})
\begin{equation}
\bL_{\vc,h}\cong \tau (\cL_{\vc}),\quad
\bL_{-\vc,h^{-1}}\cong \tau(\cL_{-\vc}).
\end{equation}
where $\tau=\mu\circ \kappa$ is the composition of
the microlocalization $\mu$ and the coherent-constructible
correspondence $\kappa$.

\subsection{T-dual Lagrangians are branes}\label{sec:compactified}
In this section, we study the behavior of T-dual Lagrangians
on the compactification  $D^*M_\bR = \overline{T}^*M_\bR$ in the cotangent.
We will show that Lagrangians $\bL_{-\vc,h}$ from anti-ample line bundles $\cL_{-\vc}$
are branes (Proposition \ref{tdualsarebranes});
as an immediate consequence, Lagrangians $\bL_{\vc, h}$ from ample line
bundles $\cL_\vc$ are also branes (Corollary \ref{thm:cset-ample}).
To prove a Lagrangian $L$ is a brane of $Fuk(T^*M_\bR;\LS),$ we need to establish
that (1) $L$ is tame, (2) $L$ has a brane structure, (3) $\pi(L)$ is bounded,
(4) $\overline{L}\subset \overline{T}^*M_\bR$ is a $\cC$-set,
(5) $L^\infty\subset \Lambda_\Si^\infty$.

\begin{proposition}[T-dual Lagrangians are tame]
\label{thm:tame}
Let $\bL_{\vc,h}$ be the T-dual Lagrangian constructed in \ref{sec:construct}.
 (We do not assume $\cL_\vc$ is ample or anti-ample.) Then:
\begin{enumerate}
\item there exists $\rho>0$ such that for every $p\in \bL_{\vc, h^{-1}}$, the
set of points $p' \in \bL_{\vc, h^{-1}}$ with $d(p, p')<\rho$ is
contractible;
\item there exists a constant $C=C(\vc,h)$ such that
$$
d_{\bL_{\vc, h^{-1}}}(p, p') < Cd(p, p')
$$
for all $p, p' \in \bL_{\vc, h^{-1}}$, where $d$ is the
distance in $T^*M_\bR$ and $d_{\bL_{\vc, h}}$ is the distance
in $\bL_{\vc, h}$.
\end{enumerate}
Therefore $\bL_{\vc,h}$ is tame in the sense of \cite{NZ}.
\end{proposition}
\begin{proof}

The Lagrangian $\bL_{\vc, h}$ is the graph of the map $\Phi_{\vc,h}\circ j_0: N_\bR \to M_\bR$. We first show that the first and second derivatives of $\Phi_{\vc,h} \circ j_0$, i.e. $\frac{\partial^2 f_{\vc,h}}{\partial y_i\partial y_j}$
and $\frac{\partial^3 f_{\vc,h} }{\partial y_i\partial y_j\partial y_l}$ are bounded for any $i,j,l$.

For each top dimensional cone $C_k\in \Sigma$, $k=1,\dots,v$, the
associated affine toric variety $U_k\cong \bC^n$ is smooth since
$X_\Si$ is a smooth projective toric variety. The coordinates in
$U_k$ are given by
$$
z_{k,i}=s_{k,i}+\sqrt{-1} t_{k,i}=r_{k,i}\exp(\sqrt{-1}\theta_{k,i})=\exp(y_{k,i}+\sqrt{-1}\theta_{k,i}).
$$
Notice that the coordinates $y_{k,i}$ and $y_i$ differ by a linear
change of basis. Fix a compact part
$U'_k=\{|z_{k,1}|^2+\dots+|z_{k,n}|^2\le M\}\subset U_k$ such that
$X_\Si=\cup_{k=1}^v U'_k$.
The 2-form
\begin{eqnarray*}
\omega_h&=&\sum_{i,j=1}^n\frac{\partial^2 f_{\vc,h} }{\partial
y_{k,i}\partial y_{k,j}} dy_{k,i}\wedge d\theta_{k,j}
= \sum_{i,j=1}^n\frac{\partial^2 f_{\vc,h} }{r_{k,i} r_{k,j}\partial
y_{k,i}\partial y_{k,j}} r_{k,j}dr_{k,i}\wedge d\theta_{k,j}\\
\ &=& \sum_{i,j=1}^n\cos(\theta_{k,i}-\theta_{k,j})\cdot
\frac{\partial^2 f_{\vc,h} }{r_{k,i} r_{k,j}\partial y_{k,i}\partial
y_{k,j}} \cdot (ds_{k,i}\wedge dt_{k,j}+ds_{k,j}\wedge
dt_{k,i})\\
\ &+& \sum_{i,j=1}^n\sin(\theta_{k,i}-\theta_{k,j})\cdot
\frac{\partial^2 f_{\vc,h} }{r_{k,i} r_{k,j}\partial y_{k,i}\partial
y_{k,j}} \cdot (ds_{k,i}\wedge ds_{k,j}+dt_{k,i}\wedge dt_{k,j}).
\end{eqnarray*}

Hence $\omega_h$ must be in the form $\omega_h= a_{k,ij} (ds_{k,i}\wedge dt_{k,j}+ds_{k,j}\wedge
dt_{k,i}) + b_{k,ij} (ds_{k,i}\wedge ds_{k,j}+dt_{k,i}\wedge
dt_{k,j})$, and we know that $a_{k,ij}$ and $b_{k,ij}$ are bounded in $U_k'$
since they are real analytic functions on $U_k$.
By comparing with the expression above,
$$ \frac{\partial^2
f_{\vc,h} }{\partial y_{k,i}\partial y_{k,j}}=\frac{a_{k,ij} r_{k,i}
r_{k,j}}{\cos(\theta_{k,i}-\theta_{k,j})}=\frac{b_{k,ij} r_{k,i}
r_{k,j}}{\sin(\theta_{k,i}-\theta_{k,j})}.$$ Thus
$$\left | \frac{\partial^2
f_{\vc,h} }{\partial y_{k,i}\partial y_{k,j}}\right |\le \sqrt{2}
\max\{|a_{k,ij}|,|b_{k,ij}|\}\cdot r_{k,i}r_{k,j}.$$ The right hand
side is bounded on $U'_k$, and therefore $\frac{\partial^2
f_{\vc,h} }{\partial y_{k,i}\partial y_{k,j}}$ is bounded on $U'_k$, for
any $i,j$.

Moreover, $$\frac{\partial^3
f_{\vc,h} }{\partial y_{k,i}\partial y_{k,j}\partial y_{k,l}}=\frac{1}{\cos(\theta_{k,i}-\theta_{k,j})}\frac{\partial(a_{k,ij} r_{k,i}
r_{k,j})}{\partial y_{k,l}}=\frac{1}{\sin(\theta_{k,i}-\theta_{k,j})}\frac{\partial(b_{k,ij} r_{k,i}
r_{k,j})}{\partial y_{k,l}}$$
also implies that on $U_k'$ the derivatives $\frac{\partial^3
f_{\vc,h} }{\partial y_{k,i}\partial y_{k,j}\partial y_{k,l}}$ are bounded since $\frac{\partial(a_{k,ij} r_{k,i}
r_{k,j})}{\partial y_{k,l}}$ and $\frac{\partial(b_{k,ij} r_{k,i}
r_{k,j})}{\partial y_{k,l}}$ are bounded on $U'_k$.

There exists constants $(C^k_{ij})$, $k=1,\ldots,v$,
such that
$$ \frac{\partial^2 f_{\vc,h}}{\partial
y_i \partial y_j} =\sum_{a,b} C^k_{ia}C^k_{jb} \frac{\partial ^2 f_{\vc,h}}
{\partial y_{k,a}\partial y_{k,b}};\ \frac{\partial^3 f_{\vc,h}}{\partial
y_i \partial y_j \partial y_l} =\sum_{a,b,c} C^k_{ia}C^k_{jb}C^k_{lc} \frac{\partial ^3 f_{\vc,h}}
{\partial y_{k,a}\partial y_{k,b} \partial y_{k,c}}.
$$
Hence there is $M_k$ such that $\displaystyle{ \Bigl|\frac{\partial^2
f_h}{\partial y_i\partial y_j}\Bigr|} <M_k$ and $\displaystyle{ \Bigl|\frac{\partial^3
f_h}{\partial y_i\partial y_j\partial y_l}\Bigr|} <N_k$ on $U_k'$ for any $i,j$. By
construction $\cup_{k=1}^v U_k'= X_\Sigma$.  It follows that for
$M=\max M_k$ and $N=\max N_k$, we have the inequalities
$$
\Bigl|\frac{\partial^2 f_{\vc,h}}{\partial y_i\partial y_j} \Bigr|< M;\ \Bigl|\frac{\partial^3 f_{\vc,h}}{\partial y_i\partial y_j \partial y_l} \Bigr|< N,
$$
for any $i,j$.

To show (1), let $p=(x_0,y_0)$ be any point in $\bL_{\vc,h}$. Let $\xi=(\xi_1,\dots,\xi_n)$ be a unit vector in $N_\bR$, and $y_t=y_0+t\xi$. Set $p_t=(x_t,y_t)\in \bL_{\vc,h}$ where $x_t=\Phi_{\vc,h}\circ j_0(y_t)$. Near $p$ the Taylor theorem gives $$x_t=x_0+ tA + t^2B(t'),$$
where $A,B$ are in $M_\bR$ with each component $$A_i=\sum_{j=1}^n \xi_j \frac{\partial^2 f}{\partial y_i \partial y_j}(y_0),\ B_i=\sum_{j,l=1}^n \xi_j \xi_l \frac{\partial^3 f}{\partial y_i \partial y_j \partial y_l}(x_{t'}),$$ and $t'\in [0,t]$ depends on $t$.
Therefore, $d(p,p_t)^2=t^2+(tA+t^2B(t'))^2$. Since by our estimates $|A|< nM$ and $|B|<n^2N^2$, there exists an $\rho>0$ such that for any direction $\xi$, $d(p,p_t)$ increases as long as $0<t<\rho$. Hence the set $\{p'\in \bL_{\vc,h}:d(p,p')<\rho\}\subset \{p'\in \bL_{\vc,h}:d_{N_\bR}(p,p')<\rho\}$ is a star-set, and it is contractible.

For any $p_1=(x_1,y_1), p_2=(x_2,y_2) \in\bL_{\vc,h}$,
$$d_{\bL_{\vc, h}}(p_1,p_2)\le \int_{l_{y_1,y_2}}
\sqrt{1+n M^2}d\xi = \sqrt{1+nM^2} d_{N_\bR}(y_1,y_2)\le
\sqrt{1+n M^2} d(p_1,p_2),$$ where $d\xi$ is the standard measure on
the segment $l_{y_1,y_2}$ from $y_1$ to $y_2$ in $N_\bR$. This shows (2).
\end{proof}

\begin{remark}
In \cite{NZ}, a new metric $g_{con}$, which is the metric of a cone over the spherical bundle $S^*M_\bR$ near the infinity, is introduced in order to ensure a tame perturbation for any standard Lagrangian. It is no longer needed here since our T-dual Lagrangians are already tame in the usual Sasaki metric. Moreover, we only consider standard or costandard Lagrangians over convex polytopes, which are also tame in the Sasaki metric. Any convex polytope is prescribed by a collection of linear functions $f_i\ge 0$ for $i=1,\dots,k$. The standard Lagrangian over it can be written as the graph of $d\log m_1 +\dots+d\log m_k$, where $m_i$ is a piecewise linear function on $M_\bR$ which is $f_i$ on the half plane $\{f_i\ge 0\}$ and zero otherwise. The tameness of this standard Lagrangian follows from the tameness of each $d\log m_i$.
\end{remark}

From now on, we assume that $\cL_\vc$ is an equivariant ample line bundle
and $\omega_h$ is symplectic.

\begin{lemma}[Compact horizontal support and brane structure] \label{brane}
$\bL_{-\vc,h^{-1}}$ and $\bL_{\vc,h}$ are horizontally compact Lagrangians inside $T^*M_\bR,$
and have canonical brane structures.
\end{lemma}
\begin{proof}
Horizontal compactness is immediate, as $\triangle_{\vc}$ and $\triangle_{-\vc}$ are bounded.
Recall that a brane structure is a relative pin structure and a choice of grading (see \cite{S3} as quoted in \cite{NZ}).
Since $\bL_{-\vc,h^{-1}}$ is the graph of
a differential $df^*_{-\vc,h^{-1}}$, for $f^*_{-\vc,h^{-1}}:\tri^\circ_{-\vc}\rightarrow \bR$
(see Section \ref{sec:ample-inverse}),
it is Hamiltonian isotopic to the zero section.\footnote{The isotopy is achieved by
he Hamiltonian flow of the function $H=
f^*_{-\vc,h^{-1}} \circ \pi$, which takes $\bL_{-\vc,h^{-1}}$ to
$(1-t) \bL_{-\vc,h^{-1}}$ in time $t$.  The subset $\tri^\circ_{-\vc}$
of the zero section is the image of time-one flow.}
Since $\tri^\circ_{-\vc}\subset 0_{T^*M_\bR}$ is a contractible subset
of the zero section, it has trivial pin structure and can be given the zero grading.
The same goes for $\bL_{\vc,h}$.
\end{proof}

We use the notation of Section \ref{sec:toric-geometry}. Given
a cone $\tau\in \Si$, define $\cU_{\tau,\pm\vc}=\pm \Phi_{\vc,h}(O_\tau^+) =\pm \Phi_{\vc,h}(X)\subset
\tri_{\pm \vc}$, where $O_\tau^+$ is defined in Section  \ref{sec:toric-geometry},
and define $F_{\tau,\pm}$ to be the closures
of $\cU_{\tau,\pm\vc}$ in $M_\bR$. Then
\begin{eqnarray*}
\cU_{\tau,\vc}&=& \{  m\in \tri_{\vc}\mid
\langle m,v_i\rangle = - c_i \Leftrightarrow v_i\in \tau\} \\
&=& \{m\in M_\bR\mid
\langle m,v_i \rangle = - c_i \ (\textup{resp.}>- c_i)
 \textup{ if } v_i \in \tau \ (\textup{resp.}\notin \tau) \} \\
F_{\tau,\vc}&=& \{ m \in \tri_{\vc}\mid
\langle m,v_i\rangle = - c_i \textup{ if } v_i\in \tau\}\\
&=&  \{m\in M_\bR\mid
\langle m,v_i \rangle = - c_i \ (\textup{resp.}\geq  -c_i)
 \textup{ if } v_i \in \tau \ (\textup{resp.}\notin \tau) \}
\end{eqnarray*}
In particular, $\cU_{\tau,\pm \vc}$ are contractible open subsets
of an affine subspace of $M_\bR$, and
$$
\cU_{\{0\},\pm \vc}= \tri_{\pm\vc}^\circ,\quad
F_{\{0\},\pm \vc}=\tri_{\pm \vc}.
$$
We have a stratification
$$
\tri_{\pm \vc}=\bigcup_{\tau\in \Si} \cU_{\tau,\pm \vc}.
$$

Given a $d$-dimensional cone $\tau\in \Si$,
$F_{\tau,\pm \vc}$ is an $(n-d)$-dimensional face of the
convex polytope $\tri_{\pm \vc}\subset M_\bR$,
and has the further stratification
$$
F_{\tau,\pm \vc} =\bigcup_{\tau\subset \sigma} \cU_{\sigma,\pm \vc}.
$$
Let $N_\tau$ be the rank $d$ sublattice of $N$ generated
by $\tau\cap N$, and let $(N_\tau)_\bR =N_\tau\otimes \bR
\cong \bR^d$.
Let $w_1,\ldots,w_d$ be defined as in Section \ref{sec:toric-geometry}, so that
$$
\tau=\Bigl\{ \sum_{j=1}^d r_j w_j \mid r_j\geq 0\Bigr\},\quad
(N_\tau)_\bR=\Bigl\{ \sum_{j=1}^d r_j w_j \mid r_j\in \bR\Bigr\}.
$$

The conormal bundle of $\cU_{\tau,\pm\vc} \subset M_\bR$ is
$$
T^*_{\cU_{\tau,\pm\vc} } M_\bR
= \cU_{\tau,\pm \vc}\times (N_\tau)_\bR  \subset M_\bR \times N_\bR
=T^* M_\bR.
$$
Its closure is the conormal bundle of $F_{\tau,\pm \vc}$\,:
$$
T^*_{F_{\tau,\pm\vc } } M_\bR  = F_{\tau,\pm \vc}\times (N_\tau)_\bR.
$$
Let $\Si'=\cup_{d>0}\Sigma(d)$, so that $\Sigma = \{ \{ 0\} \}\cup \Sigma'$.
Define  a conical Lagrangian $\Lambda_{\pm \vc}\subset T^* M_\bR$ by
$$
\Lambda_{\pm\vc} := \cU_{\{0\},\pm\vc }\times\{0\}\cup
 \bigcup_{\tau\in \Si'} \cU_{\tau,\pm\vc}\times (-\tau^\circ)
 =\bigcup_{\tau\in \Si} F_{\tau,\pm \vc}\times (-\tau)
$$
Each  $F_{\tau,\pm \vc}\times (-\tau)$ is a closed subanalytic subset of $T^*M_\bR$. Note that
$$
\Lambda_{\pm \vc} \subset \LS.
$$

Let $\iota: T^*M_\bR \to D^*M_\bR$ be defined as in \eqref{eqn:iota}.
Define
$$
\bL_{\pm \vc,h^{\pm 1}}^\infty :=\overline{\iota(\bL_{\pm \vc,h^{\pm 1 }}  )}\cap T^\infty M_\bR,
\quad
\Lambda_{\pm \vc}^\infty := \overline{\iota(\Lambda_{\pm \vc})}\cap T^\infty M_\bR.
$$
Then
$$
\Lambda_{\pm \vc}^\infty
=\bigcup_{\tau\in \Si'}\cU_{\tau,\pm \vc}\times \left( (-\tau^\circ)\cap S(N_\bR)\right)
=\bigcup_{\tau\in \Si'} F_{\tau,\pm\vc}\times \left((-\tau)\cap S(N_\bR)\right)\
$$
where $S(N_\bR)=\{ y\in N_\bR\mid |y|_{N_\bR}=1\}\cong S^{n-1}$.

We now introduce an analytic-geometric category.
(See Section  \ref{sec:review-Cset} for a brief review of
analytic-geometric categories.)
\begin{definition} \label{f-C}
Define $f:\bR\to (-1,1)$ by
\begin{equation}\label{eqn:f}
f(t)=\begin{cases}
e^{-1/t} & t>0, \\
0 & t=0, \\
-e^{1/t} & t<0.
\end{cases}
\end{equation}
Let $\cC$ be the smallest analytic-geometric category
such that $f$ is a $\cC$-map.
\end{definition}

\begin{remark}
Let $f$ be defined by \eqref{eqn:f}. Then
$f$ is $C^\infty$ on $\bR$, is
real analytic on $\bR \setminus \{0\}$,
and is a homeomorphism from $\bR$ to $(-1,1)$.
So  $f^{-1}:(-1,1)\to \bR$ is a $\cC$-map, and
$f$ is an $\cC$-isomorphism.
\end{remark}

\begin{proposition}
\label{thm:cset}
$\bL_{-\vc,h^{-1}}^\infty =\Lambda_{-\vc}^\infty$, and
$\overline{\iota(\bL_{-\vc,h^{-1}})}\subset D^*M_\bR = \overline{T}^*M_\bR$
is a $\cC$-set, where $\cC$ is the analytic-geometric category defined
in Definition \ref{f-C}.
\end{proposition}

\begin{proof}
The proof is given in Section \ref{sec:proof-Cset}.
\end{proof}

\begin{corollary}[T-dual Lagrangians are branes]
\label{tdualsarebranes} T-dual Lagrangians from anti-ample
equivariant line bundles are branes.  That is, $\bL_{-\vc,h^{-1}}$
defines an object of $Fuk(T^*M_\bR;\LS).$
\end{corollary}
\begin{proof}
First, we put the trivial vector bundle on $\bL_{-\vc,h^{-1}}$. The
existence of tame perturbations follows from Proposition
\ref{thm:tame} since one may choose the constant ``perturbation".
The remaining conditions on branes are assured by Lemma
\ref{brane} and \ref{thm:cset}.
\end{proof}

Since the involution $\beta : M_\bR\times N_\bR\to M_\bR\times N_\bR$
given by $(x,y)\mapsto (-x,y)$ is a $\cC$-isomorphism such that
$\beta(\LS)=\LS$, and the tameness
is obviously preserved, we have the immediate corollary:
\begin{corollary}
\label{thm:cset-ample}
T-dual Lagrangians from ample line bundles are branes:
$\bL_{\vc,h}$ defines an object of $Fuk(T^*M_\bR;\LS).$
\end{corollary}

\subsection{T-dual Lagrangians of ample bundles are costandard branes}
\label{sec:tlagstd}

Having shown $\bL_{-\vc,h^{-1}}$ is a brane, we now relate it
to the {\em standard} brane associated to $\tri^\circ_{-\vc}$.
The key is to study normalized geodesic flow at infinity, which
controls the hom spaces of Lagrangians which intersect at infinity.
The symplectomorphism of inversion on the fibers intertwines with Verdier duality
of constructible sheaves under microlocalization \cite{N}.  We use this fact
to relate $\bL_{\vc,h}$ to the {\em costandard} brane on the set
$\tri^\circ_\vc$.
%

\subsubsection{Normalized geodesic flow}
Let $\{e_i^*\},$ $\{e_j\}$ be dual orthonormal bases on $M$ and $N$, respectively (as in Sec. \ref{unwrapped}),
and let $x_i,$ $y_j$ be associated real coordinates.  We can equate
$e_j$ with $dx^j,$ so $(x,y) = (\sum_i x_i e^*_i, \sum_j y_j dx_j)\in M_\bR \times N_\bR = T^*M_\bR$.
The inner product on $N_\bR$ induces a linear isomorphism
$I:N_\bR\to M_\bR$ given by
$y \mapsto \sum_{j=1}^n \langle e_j^*, y\rangle e_j^*$.
In particular, $I(e_i)= e_i^*$, so $I$ is an isometry.
Define $y^*=I(y)$.
%

Given a vector space $V$, let $V' = V\setminus \{0\}$;
given a vector bundle $E$, let $E'$ denote the complement
of the zero section. The normalized geodesic flow
on $(T^*M_\bR)' \cong (TM_\bR)'$ is given by
\begin{eqnarray*}
&& \gamma_t : (T^*M_\bR)' \cong M_\bR \times N_\bR'
\to (T^*M_\bR)'\cong M_\bR\times N_\bR' ,\\
&& \gamma_t(x,y)= (x + \frac{t y^*}{|y^*|_{{M_\bR}} } , y ),
\end{eqnarray*}
where $|y^*|_{M_\bR} = |y|_{N_\bR}$ because $y\mapsto y^*$ is an isometry
from $N_\bR$ to $M_\bR$.

Let $\cL_\vc$ be an ample line bundle,
and let $h$, $\Phi_h$, $\tri_{-\vc}$, $\bL_{-\vc,h^{-1}}$,
etc.  be defined as in Section \ref{sec:ample-inverse}.
Let $q \in \pa\tri_{-\vc}\subset M_\bR$ be a boundary point of the polytope.
We consider the following two Lagrangians in $T^* M_\bR\cong M_\bR \times N_\bR$:
\begin{eqnarray*}
L_q &=&\{ (q,y)\mid y\in N_\bR\}\subset M_\bR\times N_\bR\\
\bL_{-\vc,h^{-1}}&=&\{ (-\Phi_h\circ j_0(y),y)\mid y\in N_\bR \} \subset M_\bR \times N_\bR
\end{eqnarray*}
Let $L_q' = (T^*M_\bR)'\cap L_q $, and let
$\bL_{-\vc,h^{-1}}' = (T^* M_\bR)' \cap  \bL_{-\vc,h^{-1}}$.
Then
\begin{eqnarray*}
\gamma_t( L_q' )&=&
\{ (q+\frac{t y^*}{|y^*|_{M_\bR}}, y)\mid y\in N_\bR'\}\subset M_\bR\times N_\bR' \\
\gamma_t( \bL_{-\vc,h^{-1}}' )
&=& \{ (-\Phi_h \circ j_0(y)+\frac{ty^*}{|y^*|_{M_\bR}} ,y )\mid
y\in N_\bR' \} \subset M_\bR\times N_\bR'.
\end{eqnarray*}
Note that $(x,y)\in  \gamma_{t_1}(L_q)\cap \gamma_{t_2}(\bL_{-\vc,h^{-1}} )$
if and only if
\begin{equation}\label{eqn:intersection}
y\in N_\bR',\quad  q + \Phi_h\circ j_0(y) =\frac{(t_2-t_1)y^*}{|y^*|_{M_\bR}}.
\end{equation}

\begin{lemma}
\label{nohom}
Given any $q\in \pa\tri_{-\vc}$, there
exists $\delta >0$ such that
$$
0\leq t_1\leq t_2 < \delta
\Rightarrow
\gamma_{t_1} (L_q' ) \cap
\gamma_{t_2} ( \bL_{-\vc,h^{-1}}')=\emptyset.
$$
\end{lemma}

\begin{proof} We use the notation in Section \ref{sec:compactified}.
$$
\tri_{-\vc} =\bigcup_{\tau\in \Si} \cU_{\tau,-\vc}.
$$
The right hand side is a disjoint union.
Let $\Si(d)$ be the set of $d$-dimensional cones
in $\Si$.

\noindent
{\em Step 1.} $q \in \pa\tri_{-\vc}$, so there exists a unique
$d>0$ and a unique $\tau\in \Si(d)$ such that
$q \in \cU_{\tau,-\vc}=-\Phi_h(\cO_\tau^+)$. There exists
a unique $x \in \cO_\tau^+\subset (X_\Si)_{\geq 0}$
such that $-\Phi_h(x)=q$.

There exists
$\sigma\in \Si(n)$ such that
$\tau\subset \si$. Let $w_j$ be defined
as in the proof of Proposition \ref{thm:cset} (see Section
\ref{sec:proof-Cset}), so
$$
\tau =\{ r_1 w_1 +\cdots + r_d w_d\mid r_j \geq 0\},\quad
\si= \{ r_1 w_1 +\cdots + r_n w_n\mid r_j\geq 0\}
$$
The holomorphic coordinates of
$X_\sigma =\Spec\bC[\sigma^\vee\cap M]\cong \bC^n$
are $Z_j= \chi^{w_j^\vee}$, $j=1,\ldots, n$.
There exist $b_{d+1},\ldots, b_n\in \bR$ such that
the coordinates of $x \in U_\sigma$ are given by
$$
Z_1=\cdots = Z_d=0,\quad
Z_{d+1}= e^{b_{d+1}},\ldots, Z_n= e^{b_n}.
$$

\noindent
{\em Step 2.} For any $r>0$, define
$$
S_r =\{ r_1 w_1 +\cdots+ r_n w_n \mid r_i\in (-r,r) \},\quad
B_r =\{ y\in N_\bR\mid |y|_{N_\bR} < r\}.
$$
There exists $c\in (0,1)$ such that for all $r>0$,
$$
B_{cr} \subset S_r \subset B_{c^{-1} r}.
$$
Let $R=\max\{ |b_{d+1}|,\ldots, |b_n|\}+1$.
Note that $x$ is contained in (see Section \ref{sec:toric-geometry} for definitions)
$$
X_\tau^+ :=X_\tau\cap (X_\Si)_{\geq 0}\cong [0,\infty)^d\times (\bR^+)^{n-d}
$$
which is an open set in $(X_\Si)_{\geq 0}$.
A neighborhood of $x$ in $X_\tau^+$ is given by
$$
U= \{ (Z_1,\ldots, Z_n)\mid Z_1,\ldots, Z_d\in [0, e^{-2c^{-2}R} ),\quad
Z_{d+1},\ldots, Z_{n}\in (e^{-R}, e^R)\}.
$$

Recall that $j_0: N_\bR\to U_\sigma$ is given by
$\sum_{j=1}^n r_j w_j \mapsto (e^{r_1},\ldots, e^{r_n})$, so
\begin{eqnarray*}
j_0^{-1}(U)&=&\{ r_1 w_1+\cdots +r_n w_n\mid r_1,\ldots, r_d< -2c^{-2}R,\quad
r_{d+1},\ldots, r_n\in (-R,R)\}\\
&\cong& (-\infty, -2c^{-2}R)^d\times (-R,R)^{n-d}.
\end{eqnarray*}

\noindent
{\em Step 3.}
$-\Phi_h$ maps $X_\tau^+$ homeomorphically to
$-\Phi_h(X_\tau^+)$, so there exists $\delta >0$ such that
$B(q, \delta):=\{ m\in M_\bR\mid |m-q|_{M_\bR}<\delta\} \subset -\Phi_h(U)$.

\noindent
Claim: For any $y\in N_\bR'$ and $0\leq t_1\leq t_2<\delta$,
\eqref{eqn:intersection} does not hold. Therefore,
$$
\gamma_{t_1}( L_q' ) \cap
\gamma_{t_2} ( \bL_{-\vc,h^{-1}}')=\emptyset.
$$

\noindent
Case 1. $j_0(y) \notin U$. Then $ -\Phi_h\circ j_0(y) \notin B(q,\delta)$, so
$$
|q +\Phi_h\circ j_0(y)|_{M_\bR}
=|-\Phi_h\circ j_0(y)-q|_{M_\bR} \geq \delta.
$$
On the other hand
$$
\left|\frac{(t_2-t_1) y^*}{|y^*|_{M_\bR}}\right|
= t_2-t_1 < \delta.
$$
So
$$
q + \Phi_h\circ j_0(y) \neq \frac{(t_2-t_1)y^*}{|y^*|_{M_\bR}}.
$$

\noindent
Case 2. $j_0(y)\in U$. We have
$$
y= r_1 w_1+\cdots + r_n w_n,\quad r_1,\ldots, r_d < -2c^{-2}R,\quad
r_{d+1},\ldots, r_n \in (-R,R).
$$
Let $y_1= r_1 w_1+\cdots + r_d w_d$ and $y_2= r_{d+1} w_{d+1}+\cdots r_n w_n$.
Then
$$
y=y_1+y_2,\quad
y_1\in N_\bR\setminus \overline{S_{2c^{-2}R}} \subset N_\bR \setminus \overline{B_{2c^{-1}R}},\quad
y_2 \in S_R \subset B_{c^{-1}R}.
$$
Therefore,
$$
|y_1|_{N_\bR}> 2c^{-1}R> c^{-1}R > |y_2|_{N_\bR}.
$$

Let $(v_1,v_2)_{N_\bR}$ denote the inner product on $N_\bR$, so that
$$
(e_i,e_j)_{N_\bR}=\delta_{ij},\quad \langle v_1^*, v_2\rangle = (v_1,v_2)_{N_\bR}.
$$
Then
\begin{eqnarray*}
\langle \frac{(t_2-t_1)y^*}{|y^*|_{M_\bR}}, y_1\rangle
=\frac{t_2-t_1}{|y^*|_{M_\bR}} (y_1+y_2,y_1)_{N_\bR}
\end{eqnarray*}
where $t_2-t_1\geq 0$, and
\begin{eqnarray*}
(y_1+y_2,y_1)_{N_\bR}&=& |y_1|^2_{N_\bR}+ (y_2,y_1)_{N_\bR}
\geq |y_1|^2_{N_\bR}- |y_2|_{N_\bR}|y_1|_{N_\bR}\\
&=& |y_1|_{N_\bR}( |y_1|_{N_\bR}-|y_2|_{N_\bR}) >0.
\end{eqnarray*}
So
\begin{equation}\label{eqn:positive}
\langle \frac{(t_2-t_1)y^*}{|y^*|_{M_\bR}}, y_1\rangle \geq 0.
\end{equation}

On the other hand,
\begin{equation}\label{eqn:w-sum}
\langle q +\Phi_h\circ j_0(y), y_1\rangle
=\sum_{j=1}^d r_j \langle q  + \Phi_h\circ j_0(y),w_j\rangle.
\end{equation}
Let $w_j = v_{i(j)}$. Since $q\in \cU_{\tau,-\vc}$ and $-\Phi_h\circ j_0(y) \in \tri_{-\vc}^\circ$,
for $j=1,\ldots,d$,
$$
\langle q, w_j \rangle =c_{i(j)},\quad
\langle -\Phi_h\circ j_0(y), w_j\rangle < c_{i(j)}.
$$
So we have
\begin{equation}\label{eqn:pm}
\langle q + \Phi_h\circ j_0(y),w_j\rangle >0,\quad
r_j<-2c^{-1}R< 0,
\end{equation}
Equations \eqref{eqn:w-sum} and \eqref{eqn:pm} imply
\begin{equation}\label{eqn:negative}
\langle q +\Phi_h\circ j_0(y), y_1\rangle <0.
\end{equation}

Combining \eqref{eqn:positive} and \eqref{eqn:negative}, we see that
$$
q + \Phi_h\circ j_0(y) \neq \frac{(t_2-t_1)y^*}{|y^*|_{M_\bR}}.
$$
\end{proof}

\subsubsection{Lagrangians from anti-ample line bundles are standard
branes} \label{eqn:neg-ample}


%


We now show
$\bL_{-\vc,h^{-1}}$ is isomorphic to the standard Lagrangian brane over $\tri^\circ_{-\vc}$,
and that $\bL_{\vc,h}$ is isomorphic to the costandard brane over $\tri^\circ_{\vc}.$

\begin{theorem}
\label{mainthm}
Let $\mu: Sh_{cc}(M_\bR;\LS)\rightarrow Fuk(T^*M_\bR;\LS)$ be the microlocalization quasi-embedding of
Theorem \ref{nadler-zaslow}.  Then
$\bL_{-\vc,h^{-1}} \cong \mu(i_*\bC_{\tri^\circ_{-\vc}})$,
and $\bL_{\vc, h}\cong \mu(i_! \omega_{\tri^\circ_\vc }).$
\end{theorem}

\begin{proof}
We show  $\bL_{-\vc,h^{-1}} \cong \mu(i_*\bC_{\tri^\circ_{-\vc}})$
by proving, following \cite{N}, that the two objects define isomorphic modules
under the Yoneda embedding
$$
\cY: DFuk(T^*M_\bR)\rightarrow mod(DFuk(T^*M_\bR)),\qquad
\cY(L) = hom_{DFuk(T^*M_\bR)}(-,L).$$
To prove that $\cY(\bL_{-\vc,h^{-1}})\cong \cY( \mu(i_*\bC_{\Delta^\circ_{-\vc}})),$
we first fix a triangulation $\mathcal T$ of $M_\bR$
containing $\{\cU_{\tau,-\vc}\mid \tau\in\Sigma\}$
(recall $\cU_{\{0\},-\vc}= \tri^\circ_{-\vc}$).
The technique of \cite{N} exploits the triangulation to resolve the
diagonal standard, i.e. the identity functor.
What emerges is that the Yoneda module of any object $\cY(L)$ is expressed
in terms of (sums and cones of shifts of) Yoneda modules from standards,
$\cY(\mu(i_*\bC_{T})),$ where $T\in \cT.$
The coefficient of the Yoneda standard module $\cY(\mu(i_*\bC_{T})),$
takes the form $hom_{DFuk(T^*M_\bR)}(L_{\{t\}*},L),$
where $t$ is any point in $T$ (contractibility of $T$ means
that the choice is irrelevant up to isomorphism)---see Remark 4.5.1 of \cite{N}.

We now apply this to $L=\bL_{-\vc,h^{-1}}.$
First consider the case $T\neq \tri^\circ_{-\vc},$ and let $t\in T.$
Then if $T\cap \tri_{-\vc} = \emptyset,$
clearly $hom_{DFuk(T^*M_\bR)}(L_{\{t\}*},\bL_{-\vc,h^{-1}}) = 0,$
since $L_{\{t\}*}$ is just the fiber $T^*_tM_\bR.$
Otherwise, if
$T\cap \partial \tri$ is nonempty, then Proposition \ref{nohom}
ensures us that $hom_{DFuk(T^*M_\bR)}(L_{\{t\}*},\bL_{-\vc,h^{-1}}) = 0.$
Finally, if $T = \tri^\circ_{-\vc},$ then since $\bL_{-\vc,h^{-1}}$ is
a graph over $T,$ we have
$hom_{DFuk(T^*M_\bR)}(L_{\{t\}*},\bL_{-\vc,h^{-1}}) = \bC.$
Therefore, $\bL_{-\vc,h^{-1}} \cong \mu(i_*\bC_{\tri^\circ_{-\vc}}),$
and the first statement is proved.
Note that the result is independent of how $\mathcal T$ was chosen.

The map $\alpha: (x,y)\mapsto (x,-y)$ gives rise to a duality functor (still denoted by $\alpha$)
$$
\alpha: Fuk(T^*M_\bR)^\circ \to Fuk(T^*M_\bR).
$$
The functor $\alpha$ sends a Lagrangian brane $L$ to $\alpha(L)$. It is proved in Section 5.1 of \cite{N} (Proposition 5.1.1) that there is a functor quasi-isomorphism
$$
\mu\circ \cD\cong \alpha\circ \mu: Sh_{cc}(M_\bR)\to TrFuk(T^*M_\bR).
$$

Define another functor $\nu: Fuk(T^*M_\bR)\to Fuk(T^* M_\bR)$ given by the map
$$
M_\bR\times N_\bR\to M_\bR\times N_\bR,\quad (x,y)\mapsto (-x,-y).
$$
The functor $\nu$ maps any standard brane $L(U)$ over the submanifold $U\hookrightarrow M_\bR$ to the standard brane $L(-U)$ over $-U$. Let $\cR$ be the induced push-forward on $Sh_c(M_\bR)$ given by the map $x\mapsto -x$. It is obvious that there is an isomorphism of functors:
$$
\mu\circ \cR \cong \nu \circ \mu: Sh_c(M_\bR)\to TrFuk(T^* M_\bR).
$$
Therefore, the quasi-isomorphism $\bL_{-\vc,h^{-1}}\cong \mu(i_* \bC_{\tri_{-\vc}^\circ})$ gives rise to
$$
\nu(\bL_{-\vc,h^{-1}})\cong \nu(\mu(i_* \bC_{\tri_{-\vc}^\circ})) \cong
\mu(\cR(i_* \bC_{\tri_{-\vc}^\circ}))\cong \mu(i_* \bC_{\tri_{\vc}^\circ}).
$$
The quasi-isomorphism $\mu\circ \cD\cong \alpha\circ \mu$
induces
$$
\alpha(\nu(\bL_{-\vc,h^{-1}}))\cong\alpha(\mu(i_* \bC_{\tri_{\vc}^\circ}))
\cong \mu(\cD(i_* \bC_{\tri_{\vc}^\circ })) \cong \mu(i_! \omega_{\tri_{\vc}^\circ}).$$
It is easy to see that $\alpha(\nu(\bL_{-\vc,h^{-1}}))=\bL_{\vc, h}$. Therefore we have
$$
\bL_{\vc, h}\cong \mu(i_! \omega_{\tri_{\vc}^\circ }).
$$
\end{proof}

\begin{appendix}

\section{Review of Geometric Categories and Proof of Proposition \ref{thm:cset}}
\label{sec:appCset}

\subsection{Review of analytic-geometric categories}
\label{sec:review-Cset}

We recall definitions and basic properties from \cite{vM}.

\begin{definition}[analytic-geometric category]\label{def:Cset-app}
We say that an {\em analytic-geometric category}
$\cC$ is given if each manifold $X$ is equipped with a collection
$\cC(X)$ of subsets of $X$ such that the following conditions are satisfied
for all manifolds $X$ and $Y$:
\begin{enumerate}
\item[AG1.] $\cC(X)$ is a Boolean algebra of subsets of $X$, with $X\in \cC(M)$.
\item[AG2.] If $A\in \cC(X)$, then $A\times \bR\subset \cC(X\times \bR)$.
\item[AG3.] If $f: X\to Y$ is a proper analytic map and $A\in \cC(X)$,
then $f(A)\in \cC(Y)$.
\item[AG4.] If $A\subset X$, and $(U_i)$ is an open covering of $X$
($i$ in some index set $I$), then
$A\in \cC(X)$ if and only if $A\cap U_i\in \cC(U_i)$ for all $i\in I$.
\item[AG5.] Every bounded set in $\cC(\bR)$ has finite boundary.
\end{enumerate}
\end{definition}

It is proved in \cite[Appendix D]{vM} that
this indeed gives rise to a category $\cC$.
An object of $\cC$ is a pair $(A,X)$ with $X$ a manifold and $A\in \cC(X)$.
A morphism $(A,X)\to (B,Y)$ is a continuous map $f:A\to B$ whose graph
$$
\Gamma_f =\{ (a, f(a))\mid a\in A\} \subset A\times B
$$
belongs to $\cC(X\times Y)$. We usually refer
to an object $(A,X)$ of $\cC$ as the $\cC$-{\em set} $A$ in $X$, or even
just the $\cC$-set $A$ if its ambient manifold is clear from context.
Similarly, a morphism $f:(A,X)\to (B,Y)$ is called a
$\cC$-{\em map} $f:A\to B$ if $X$ and $Y$ are clear from context.

The following basic properties are proved in \cite[Appendix D]{vM}.
\begin{theorem}\label{thm:basics}
Let $X$, $Y$ be manifolds of dimension $m$, $n$, respectively,
and let $A\in \cC(X)$, $B\in\cC(Y)$.
\begin{enumerate}
\item Every analytic map $f: X\to Y$ is a $\cC$-map.
\item Given an open covering $(U_i)$ of $X$, a map $f:A\to Y$ is a $\cC$-map
if and only if each restriction $f|_{U_i\cap A}: U_i\cap A\to Y$ is a $\cC$-map.
\item $A\times B\in\cC(X \times Y)$, and the projections
$A\times B\to A$ and $A\times B\to B$ are $\cC$-maps.
\item If $f:A\to Y$ is a proper $\cC$-map and $Z\subset A$
is a $\cC$-set, then $f(Z)\in \cC(Y)$.
\item If $A$ is closed in $X$ and $f:A\to Y$ is a $\cC$-map,
then $f^{-1}(B)\in \cC(X)$.
\item If $B_1,\ldots,B_k$ are $\cC$-sets (in possibly different manifolds),
then a map
$$
f=(f_1,\ldots,f_k): A\to B_1\times \cdots \times B_k
$$
is a $\cC$-map if and only if each $f_i:A\to B_i$ is a $\cC$-map.
\item $cl(A), int(A) \in \cC(X)$.
\end{enumerate}
\end{theorem}

\begin{corollary}
Assume $f:X\to Y$ is a $\cC$-map which is also a homeomorphism.
Then
\begin{enumerate}
\item $f^{-1}: Y\to X$ is a $\cC$-map.
\item For any subset $A\subset X$, $A\in \cC(X)
\Leftrightarrow f(A)\in \cC(Y)$.
\end{enumerate}
\end{corollary}

\subsection{Proof of Proposition \ref{thm:cset}}
\label{sec:proof-Cset}

\begin{proof}
By (7) of Theorem \ref{thm:basics},
it suffices to prove that
\begin{enumerate}
\item[(1a)]  $\iota(\bL_{-\vc,h^{-1}})\in \cC(D^*M_\bR)$,
\item[(1b)]  $\bL_{-\vc, h^{-1}}^\infty = \Lambda^\infty_{-\vc}$.
\end{enumerate}

Let
\begin{eqnarray*}
B(N_\bR) &=& \{ y\in N_\bR\mid |y|_{N_\bR} <1 \},\\
\bar{B}(N_\bR)&=& \{ y\in N_\bR \mid |y|_{N_\bR}\leq 1\}
= B(N_\bR)\cup S(N_\bR).
\end{eqnarray*}
Define
$$
F: X_\Si \times \bar{B}(N_\bR) \to M_\bR\times \bar{B}(N_\bR) = D^*M_\bR,\quad
(x,y)\mapsto (-\Phi_h(x),y).
$$
Then $F$ is a proper real analytic map. Let
$$
L=\{ (x,y)\in X_\Sigma\times  B(N_\bR)\mid x= j_0\Bigl(\frac{y}{\sqrt{1-|y|^2_{N_\bR}} }\Bigr)\}.
$$
Let $\bar{L}$ be the closure of $L$ in $\overline{B}(N_\bR)$, and let
$L^\infty =\bar{L} \cap (X_\Si\times S(N_\bR))$.  By (4)
of Theorem \ref{thm:basics},
it suffices to prove that
\begin{enumerate}
\item[(2a)] $L\in \cC(X_\Si\times \bar{B}(N_\bR))$,
\item[(2b)] $L^\infty =\displaystyle{\bigcup_{\tau\in \Si'}
O^+_\tau\times ((-\tau^\circ)\cap S(N_\bR) ) }$.
\end{enumerate}
Recall that $X_\si \cong \bC^n$ for $\si\in \Si(n)$, and
$\{ X_\si \mid \si \in \Si(n)\}$ is an open cover of
$X_\Si$. By AG4 of Definition \ref{def:Cset-app}, it suffices to prove
that, for any $\si\in \Si(n)$,
\begin{enumerate}
\item[(3a)] $L\cap (X_\si\times \bar{B}(N_\bR)) \in
\cC(X_\si\times \bar{B}(N_\bR))$,
\item[(3b)] $L^\infty \cap (X_\si\times \bar{B}(N_\bR))
=\displaystyle{ \bigcup_{\tau\in \Si', \tau\subset \sigma}
O^+_\tau\times ((-\tau^\circ)\cap S(N_\bR))}$.
\end{enumerate}

Given $\si \in \Si(n)$, there exists a $\bZ$-basis
$\{ w_1,\ldots, w_n\}$ of $N$ such that
\begin{eqnarray*}
&& \{ w_,\ldots,w_n\} \subset \{ v_1,\ldots, v_r\},\\
&& \sigma=\{ r_1 w_1 +\cdots + r_n w_n\mid r_j\geq 0\}.
\end{eqnarray*}
Let $\{ w_1^\vee,\ldots, w_n^\vee\}$ be the dual $\bZ$-basis
of $M$, so that
$$
\si^\vee =\{ s_1 w_1^\vee +\cdots + s_n w_n^\vee \mid s_j \geq 0\}.
$$
We have
$$
\bC[\si^\vee\cap M] = \bC[\chi^{w_1^\vee},\ldots, \chi^{w_n^\vee}].
$$
Let $Z_j=\chi^{w_j^\vee}$. Then $Z_1,\ldots, Z_n$ are  holomorphic
coordinates of
$$
X_\si =\Spec\bC[\si^\vee\cap M] \cong \bC^n.
$$
The image of $j_0: N_\bR\to X_\Si$ is contained
in $X_\si$, and $j_0$ is given by
$$
y\mapsto  (e^{\langle w_1^\vee,y\rangle},\ldots, e^{\langle w_n^\vee,y\rangle}),
$$
or equivalently,
$$
\sum_{j=1}^n y_j w_j \mapsto (e^{y_1},\ldots, e^{y_n}).
$$
Let $(\ , \ )_{N_\bR}$ denote the inner product on $N_\bR$,
and let $g_{ij}= ( w_i,w_j)_{N_\bR}$. Define
$Q:\bR^n\to \bR$ by
$$
Q(y_1,\ldots,y_n) \define\Bigl|\sum_{j=1}^n y_j w_j \Bigr|^2_{N_\bR}
= \sum_{j,k=1}^n g_{jk} y_j y_k.
$$
Define
$$
\bar{B}=\{ y\in\bR^n \mid Q(y) \leq 1\},
\quad S=\{ y\in\bR^n \mid Q(y)=1\}.
$$
Then $\bar{B}$ is a solid ellipsoid in $\bR^n$. Define
$$
\psi: \bR^n\times \bar{B} \longrightarrow X_\si\times \bar{B}(N_\bR),\quad
(x,y) \mapsto
(x, \sum_{j=1}^n y_j w_j)
$$
where $x=(x_1,\ldots,x_n)$ and $y=(y_1,\ldots, y_n)$.
Then $\psi$ is an injective, proper, real analytic map. Define
\begin{eqnarray*}
L_1 := \psi^{-1}(L)
&=& \{ (x, y)\in \bR^n\times \bar{B}\mid Q(y) < 1,\  x_i = \exp\Bigl(\frac{y_i}{\sqrt{1-Q(y)} }\Bigr) \}\\
&=& \Bigl\{ (x,y)\in \bR^n \times \bar{B} \mid
Q(y) <1,\   x_j >0, \  \frac{y_i}{\log x_i} =\sqrt{1-Q(y)}  \Bigr \}.
\end{eqnarray*}
Let $\bar{L}_1$ be the closure of $L_1$ in $\bR^n\times\bar{B}$,
and let $L_1^\infty = \bar{L}_1 \cap (\bR^n\times S)$.
Then
$$
\psi(L_1)=L,\quad
\psi(\bar{L}_1)=\bar{L},\quad
\psi(L^\infty_1)=L^\infty.
$$
So it suffices to prove that
\begin{enumerate}
\item[(4a)] $L_1 \in \cC(\bR^n\times \bar{B})$,
\item[(4b)] $L_1^\infty=\{ (x,y)\in\bR^n\times S \mid
y_j\leq 0,\quad x_j \geq 0,\quad  x_1 y_1=\cdots= x_n y_n =0\}$.
\end{enumerate}

Given any subset $I\subset \{1,2,\ldots,n\}$, define
$$
U_I=\{ (x_1,\ldots,x_n) \mid |x_i|<1 \textup{ for }i\in I,\
|x_i|>\frac{1}{2} \textup{ for }i\notin I\}.
$$
Then $\{ U_I \mid I\subset \{1,\ldots,n\} \}$ is an open
cover of $\bR^n$, and
$\{ U_I\times \bar{B} \mid I\subset\{1,\ldots,n\} \}$
is an open cover of $\bR^n\times \bar{B}$.
By AG4  of Definition \ref{def:Cset-app},
it suffices to prove that, for any $I\subset \{ 1,\ldots,n\}$,
\begin{enumerate}
\item[(5a)]$L_1\cap (U_I\times \bar{B} ) \in \cC(U_I\times \bar{B})$.
\item[(5b)]
\begin{eqnarray*}
L_1^\infty \cap (U_I \times\bar{B})
&=&\{ (x,y)\in  U_I\times S \mid x_j\geq 0\textup{ for } j=1,\ldots,n \\
&&   y_i\leq 0\textup{ and } x_iy_i=0 \textup{ for } i\in I,\quad
 y_i=0 \textup{ for } i\notin I \}
\end{eqnarray*}
\end{enumerate}

Without of loss of generality, we assume that
$I=\{ 1, 2,\ldots, d\}$, where $0\leq d\leq n$.
(In particular, $I$ is empty when $d=0$.)
The other cases can be obtained by permutations of
$\{1,\ldots,n\}$.

Let $J= (-\infty,-1/2)\cup (1/2,\infty)$. Then
$U_I = (-1,1)^d \times J^{n-d}$. Define
\begin{eqnarray*}
\phi:\bR^d \times J^{n-d}\times \bar{B}
&\longrightarrow& U_I \times \bar{B} \\
((t_1,\ldots,t_d), (x_{d+1},\ldots, x_n),y)
&\mapsto& ((f(t_1),\ldots,f(t_d)), (x_{d+1},\ldots, x_n),y).
\end{eqnarray*}
Then $\phi$ is a homeomorphism, and both
$\phi$ and $\phi^{-1}$ are $\cC$-maps.
To prove (5a), it suffices to prove that
$$
\phi^{-1}(L_1\cap (U_I\times \bar{B}))
\in \cC(\bR^d\times J^{n-d}\times \bar{B}).
$$
We will prove that
$$
\phi^{-1}(L_1\cap (U_I\times \bar{B}))
\in \cC_{an}(\bR^d\times J^{n-d}\times \bar{B}).
$$

We have
\begin{eqnarray*}
&& \phi^{-1}(L_1\cap (U_I\times \bar{B}) )\\
&=&\{ ((t_1,\ldots,t_d),(x_{d+1},\ldots, x_n), y)
\in \bR^d \times J^{n-d}\times \bar{B}\mid
t_i>0, x_i >\frac{1}{2}, Q(y)<1,\\
&& \sqrt{1-Q(y)} = -t_i y_i, \ 1\leq i\leq d;\
\log x_i\sqrt{1-Q(y)} = y_i,\ d+1\leq i\leq  n \}\\
&=&\{ ((t_1,\ldots,t_d),(x_{d+1},\ldots, x_n), y)
\in \bR^d \times J^{n-d}\times \bar{B}\mid
t_i>0, x_i >\frac{1}{2},\\
&&  y_1,\ldots,y_d<0,\quad (\log x_{d+1}) y_{d+1},\ldots,
(\log x_n) y_n \geq 0,\quad   Q(y) <1\\
&& 1-Q(y) = t_i^2 y_i^2, \ 1\leq i\leq d;\quad
(\log x_i)^2 (1-Q(y)) = y_i^2,\ d+1\leq i\leq  n \}
\end{eqnarray*}
$\phi^{-1}(L_1 \cap (U_I\times\bar{B}))$ is
defined by equalities and inequalities of
real analytic functions, so
$$
\phi^{-1}(L_1 \cap (U_I\times\bar{B})) \in \cC_{an}(\bR^d\times J^{n-d}\times \bar{B}).
$$
This proves (5a).

Note that
\begin{eqnarray*}
&& \phi^{-1}(\bar{L}_1 \cap (U_I\times \bar{B}) )\\
&=&\{ ((t_1,\ldots,t_d),(x_{d+1},\ldots, x_n), y)
\in \bR^d \times J^{n-d}\times \bar{B}\mid
t_i\geq 0, x_i >\frac{1}{2},\\
&&  y_1,\ldots,y_d\leq 0,\quad (\log x_{d+1}) y_{d+1},\ldots, (\log x_n) y_n
\geq 0,\\
&& 1-Q(y)= t_i^2 y_i^2, \ 1\leq i\leq d;\quad
(\log x_i)^2 (1-Q(y)) = y_i^2,\ d+1\leq i\leq  n \}
\end{eqnarray*}
which is a $\cC_{an}$-set in $\bR^d\times J^{n,d}\times \bar{B}$.
\begin{eqnarray*}
&& \phi^{-1}(L_1^\infty\cap (U_I\times \bar{B}) )\\
&=&\{ ((t_1,\ldots,t_d),(x_{d+1},\ldots, x_n), y)
\in \bR^d \times J^{n-d}\times S \mid
t_i\geq 0, x_i >\frac{1}{2},\\
&&  y_1,\ldots,y_d\leq 0,\  t_i y_i =0 , \ 1\leq i\leq d;\quad y_i=0,\ d+1\leq i\leq  n \}.
\end{eqnarray*}
This proves (5b).
\end{proof}

\section{Generating sets of line bundles}
\label{sec:ample-generate}
The Lagrangians that generate our Fukaya category are T-dual to equivariant ample line bundles.  In this section we will show that such line bundles also generate the category of equivariant coherent sheaves.  The theorem we are after is a slight generalization of a theorem of Seidel:

\begin{theorem}[Seidel]
\label{thm:seidel}
If $X$ is smooth and projective, then $\Perf_T(X)$ is generated by line bundles.
\end{theorem}
\begin{proof}
The proof of  \cite[Proposition 1.3]{Ab2} shows that $\Perf(X)$ is generated by
line bundles.  The same proof works in the $T$-equivariant setting, by the following observation.
Given a $T$-equivariant coherent sheaf $\cF$, there exists a $T$-equivariant
ample line bundle $\cL$ such that the underlying nonequivariant coherent
sheaf $\cF\otimes \cL$ is generated by global sections. 
The $T$-action on $\cF\otimes \cL$ induces
a $T$-action on  $H^0(X,\cF\otimes \cL)$. 
There exists a basis $s_1,\ldots, s_N$ of
$H^0(X,\cF\otimes \cL)$ and characters $\chi_1,\ldots, \chi_N\in \Hom(T,\bC^*)$
such that $t\cdot s_i = \chi_i(t)s_i$ for all $t\in T$. Then $s_1,\ldots, s_N$ defines
a surjective morphism  $\oplus_{i=1}^N \cL^{-1}\otimes \cO_X(\chi_i) \to \cF$ of $T$-equivariant
coherent sheaves, where  $\cO_X(\chi_i)$ is the structure sheaf equipped with
the $T$-equivariant structure given by the character $\chi_i$.
\end{proof}

The stronger version we prove is the following:

\begin{theorem}
\label{ampgen}
If $X$ is smooth and projective, then $\Perf_T(X)$ is generated by
$T$-equivariant ample line bundles.
\end{theorem}

\begin{proof}
Let $\cA$ be the full triangulated dg sub-category of $\Perf_T(X)$ generated
by $T$-equivariant ample line bundles.  We need to show that $\cA = \Perf_T(X)$.
We may see that $\cA$ is a full, dense triangulated subcategory of $\Perf_T(X)$ by the same argument used in the proof of Theorem \ref{thm:seidel} given in \cite{Ab2}. 
 (Recall that a triangulated subcategory is called \emph{dense} if every object is a direct summand of an object in the subcategory.)

Now by \cite[Theorem 2.1]{Th}, to show that $\cA = \Perf_T(X)$ it suffices to show that
the subgroup $K(\cA)$ of $K(\Perf_T(X))=K_T(X)$ is equal to $K_T(X)$.
We will show that $K_T(X)$ is additively generated
by $T$-equivariant ample line bundles.

Let $r = |\Sigma(1)|$ be the number
of $1$-dimensional cones, and let $v = |\Sigma(n)|$ be the number
of maximal cones, which is also equal to the number
of $T$-fixed points in $X$. Then
$r = \mathrm{rank}_\bZ \Pic_T(X)$.

\paragraph{\em Step 1.} Claim: There exists a
$\bZ$-basis $\{ L_1,\ldots, L_r\}$ of
$ (\Pic_T(X), \otimes)$ such that
$L_1,\ldots, L_r$ are $T$-equivariant ample line bundles.

There exists a primitive ample class $\alpha\in H^{1,1}(X;\bZ)$.
Let $M_1$ be a $T$-equivariant line bundle with
$c_1(M_1)=\alpha$. There exist $T$-equivariant line bundles
$M_2,\ldots, M_r$ such that
$\{ M_1,\ldots, M_r \}$ is a $\bZ$-basis of
$( \Pic_T(X), \otimes)$. There exist positive integers
$n_2,\ldots, n_r$ such that
$M_i \otimes M_1^{\otimes n_i}$ are ample, $i=2,\ldots,r$. Let
$$
L_1 = M_1, \quad L_i = M_i\otimes M_1^{\otimes n_i} \textup{ for } i=2,\ldots, r.
$$
Then $\{ L_1,\ldots, L_r\}$ is the desired $\bZ$-basis
of $(\Pic_T(X),\otimes)$.

\paragraph{\em Step 2.}  Let $e_i = (c_1)_T(L_i) \in H^2_T(X;\bZ)$.
Let $x_1,\ldots, x_v$ be the $T$-fixed points of $X$, and let
$\epsilon_j : x_j \to X$ be the inclusion. Let
$$
u_{ij} = \epsilon_j^* e_i \in H^2_T (x_j;\bZ) \cong M.
$$
Let $r_j: X\to x_j$ be the constant map. This gives
rise to $r_j^*: H^2_T (x_j;\bZ) \cong M \to H^2_T(X;\bZ)$.
Therefore we may view $u_{ij}\in M$ as elements in
$H_T^2(X;\bZ)$. The map $\epsilon_j^* \circ r_j^* : M\to M$ is the
identity map.  By localization, for $i=1,\ldots,r$, we have
the following relation in $H^*_T(X)$:
$$
\prod_{j=1}^v (e_i - u_{ij})=0.
$$

Let $V_{ij}$ be the $T$-equivariant line
bundle with $(c_1)_T(V_{ij}) = - u_{ij}$.
It is a $T$-equivariant lifting of the
trivial holomorphic line bundle $\cO_X$.
Define
$$
y_{ij} = \mathrm{ch}_T(L_i\otimes V_{ij})
=e^{e_i - u_{ij}},\quad
i=1,\ldots,r,\quad j=1,\ldots,v.
$$
Then we have
$$
\prod_{j=1}^v(y_{ij}-1)=0 \in H^*_T(X;\bQ),\quad i=1,\ldots,r,
$$
so
\begin{equation} \label{eqn:LV}
\prod_{j=1}^v(L_i \otimes V_{ij}-1) =0 \in K_T(X), \quad i=1,\ldots,r.
\end{equation}

\paragraph{\em Step 3.} By \cite[Proposition 3]{M}, any element in
$K_T(X)$ can be written as
\begin{equation}\label{eqn:a}
\sum a_{m_1,\ldots, m_r} L_1^{\otimes m_1}\otimes \cdots \otimes L_r^{\otimes m_r},
\end{equation}
where
\begin{enumerate}
\item[(i)] $m_1,\ldots, m_r, a_{m_1,\ldots, m_r}$ are integers, and
\item[(ii)] all but finitely many $a_{m_1,\ldots,m_r}$ are zero.
\end{enumerate}

We may use \eqref{eqn:LV} to rewrite \eqref{eqn:a} as
\begin{equation}
\sum b_{m_1,\ldots, m_r} L_1^{\otimes m_1}\otimes\cdots\otimes L_r^{\otimes m_r},
\end{equation}
where
\begin{enumerate}
\item[(i)'] $m_1,\ldots, m_r \in \{ 1, 2,\ldots, v\}$ (in particular,
the sum is finite), and
\item[(ii)'] $b_{m_1,\ldots, m_r}\in \bZ[M]$, the representation ring of $T$.
\end{enumerate}
Note that (i)' implies that, for any equivariant lifting $V$ of the trivial
holomorphic line bundle $\cO_X$, $V\otimes L_1^{ \otimes m_1}\otimes \cdots\otimes L_r^{\otimes m_r}$
is ample. Therefore $K_T(X)$ is additively generated by
$T$-equivariant ample line bundles.
\end{proof}

\section{Relation to Other Work}
\label{sec:appothers}

\subsection{Seidel, Auroux-Katzarkov-Orlov}
\label{sec:seidelako}



The homological mirror symmetry proofs of Seidel \cite{S2} and
Auroux-Katzarkov-Orlov \cite{AKO1,AKO2}, formulated in the
Fukaya-Seidel version of the mirror, make use of the fact that the
mirror categories are generated by a finite collection of objects
(Lagrangian thimbles). Studying the images of a generating set (such
as $\cO(-1),$ $\cO,$ $\cO(1)$ for $\bP^2)$) in different
formulations of homological mirror symmetry leads to the conjecture
that the thimbles are equivalent as objects to the T-dual branes
associated to these line bundles. More generally, one should search
for a proof that the dictionary between superpotential $W_\Sigma$
(see Section \ref{sec:abouzaid}) and microlocal condition
$\Lambda_\Si$ leads to an equivalence of categories.

{\bf Example: The projective plane $\bP^2$}. The mirror
Landau-Ginzburg model of $\bP^2$ is $(\bC^*)^2$ together with the
superpotential $W=z_1+z_2+1/{z_1z_2}$. The three critical points are
$(1,1),(w,w),(\bar w,\bar w)$ where
$w=-\frac{1}{2}+\frac{\sqrt{-3}}{2}$, over the critical values $3,
3w, 3\bar w$ respectively. The Fukaya-Seidel category is a category
of Lagrangian thimbles $T_i$ together with directed perturbation
when computing morphisms. These infinite Lagrangian branes $T_i$ are
the clockwise labeled vanishing thimbles over the positive-pointing
rays $\lambda_i$ starting from the critical values of $W$, parallel
to the real axis.
The difference between T-dual costandard/standard Lagrangians branes in
$Fuk(T^*T^\vee_\bR;\bar{\Lambda}_\Sigma)$ and Lagrangian thimbles in
$FS((\bC^*)^n,W)$ is illustrated in the following figure by looking
at their images under the superpotential $W$. Indeed, the T-dual
branes are very much like the vanishing thimbles in the case of
$\bP^2$: their images under $W$ also propagate from the critical
values, but in a ``thickened" way.

\smallskip
\begin{center}
\psfrag{O}{\Tiny $\cO$} \psfrag{O(1)}{\Tiny $\cO(1)$}
\psfrag{O(-1)}{\Tiny $\cO(-1)$}
\includegraphics[scale=0.35]{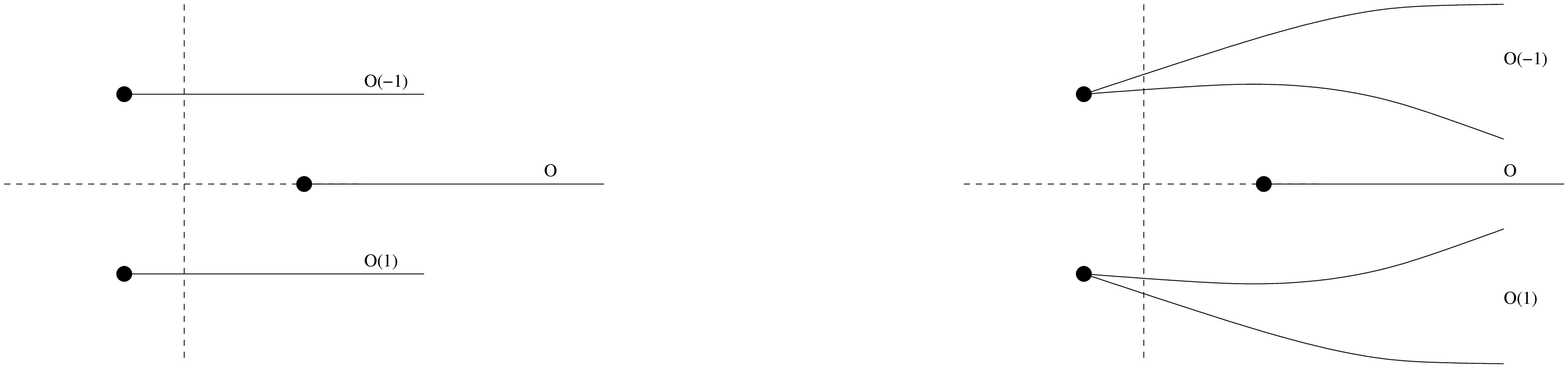}
\parbox{10cm}{Fig.2 The $W$-plane of Lagrangian A-branes in the
mirror Landau-Ginzburg model of $\bP^2$. The images of the
Lagrangian thimbles under the superpotential $W$ are horizontal rays
toward positive infinity, shown on the left. The images of the
T-dual Lagrangians (with respect to $\cO(-1),\cO,\cO(1)$) are shown
on the right, which are the areas inside the curves. They are
``thickened" versions of Lagrangian thimbles. Dashed lines are
coordinate axis.}
\end{center}


\subsection{Abouzaid}
\label{sec:abouzaid}

Abouzaid studies the Fukaya category of the Landau-Ginzburg
model dual to the toric variety.
We will describe the construction in \cite{Ab2}
in our notation. We use the notation of Section \ref{sec:toric-geometry}.
$X_\Si$ is an $n$-dimensional smooth projective toric
variety defined by a smooth complete fan $\Sigma \subset N_\bR$,
$\Si(1)=\{ \rho_1,\ldots, \rho_r\}$ is the set of 1-dimensional
cones in $\Si$, and $v_i\in N$ is the generator of
$\rho_i$, $i=1,\ldots,r$.

Let $P\subset N_\bR$ be the convex hull of $\{ v_1,\ldots,v_r\}$. The Landau-Ginzburg
model dual to $X_\Si$ is a pair $( (\bC^*)^n ,W)$, where
$W:(\bC^*)^n\to \bC$ is known as the superpotential. In our notation,
$W$ is a holomorphic function on $T^\vee = M\otimes \bC^*$,
the complex dual to the torus $T= N\otimes \bC^*$ acting
on $X_\Si$. Let $z^\alpha \in \Hom(T^\vee,\bC^*)$ be the image
of $\alpha\in N$ under the isomorphism $N \stackrel{{}_\sim}{\to} \Hom(T^\vee,\bC^*)$.
The superpotential $W:T^\vee \to \bC$ is a Laurent polynomial
$$
W =\sum_{\alpha \in N} c_\alpha z^\alpha,\quad c_\alpha \in \bC
$$
with the constraint
\begin{equation}\label{eqn:WP}
\mathrm{Newton}(W) := \{ \alpha \mid c_\alpha \neq 0\} = P.
\end{equation}

Up to now, $W$ depended only on the fan  $\Si$.  To apply
tropical geometry, Abouzaid picks an ample line
bundle $\cL_\nu$ on $X_\Si$ associated to
a strictly convex piecewise linear function $\nu:N_\bR \to \bR$
and defines a 1-parameter family of superpotentials
(recall that $\nu$ takes integral values on the lattice $N$):
$$
W_t = \sum_{\alpha \in N} c_\alpha  t^{-\nu(\alpha)} z^\alpha,\quad t\in \bC^*,
$$
where $\{ c_\alpha \}$ are fixed constants satisfying \eqref{eqn:WP}.
Therefore $( (\bC^*)^n, W_t)$ can be viewed as the dual
of the polarized toric variety $(X_\Si, \cL_\nu)$.
$M_t = W_t^{-1}(0)$ is a smooth hypersurface in $(\bC^*)^n$.

We have
$$
T^\vee \cong M\otimes \bC^* \cong  (M_\bR/M) \times M_\bR \cong T(T^\vee_\bR)
$$
Under the isomorphism $T^\vee \cong \bC^*$, the projection
$T^\vee_\bR \times M_\bR  \to M_\bR$  gets identified with
the logarithm map $\mathrm{Log}: (\bC^*)^n \to \bR^n$ in tropical
geometry:
$$
\mathrm{Log}(z_1,\ldots, z_n)  = (\log|z_1|,\ldots, \log|z_n|).
$$
Let $\cA_t := \mathrm{Log}(M_t)$ be the amoeba of $M_t$.
When $X_\Si$ is Fano, there is a unique {\sl bounded}
connected component $Q_t$ of $\bR^n - \cA_t$.  Abouzaid
defines a pre-category of {\em tropical Lagrangian sections}
whose objects
(Lagrangian branes) are sections of the restrictions of the
logarithm moment map to $Q_t$;  these Lagrangian
branes are compact $n$-dimensional
submanifolds of $T_\bR^\vee \times M_\bR$ with boundary in $M_t$.

We will describe the relation between the tropical version
of Abouzaid's Lagrangian branes (see  \cite[Section 3.3]{Ab1}) and ours.
Let
$$
\Pi = \lim_{t\to \infty} \frac{ \mathcal{A}_t}{\log t} \subset  M_\bR \cong \bR^n
$$
be the tropical amoeba, let $Q\subset M_\bR$ and  $M_\infty \subset
T(M_\bR/M)\cong (\bC^*)^n$ be the corresponding
limits of $Q_t$ and $M_t$ as $t\to \infty$. Then $Q$ is a connected
component of $\bR^n -\Pi$. Indeed, $Q$ is the moment
polytope of the ample line bundle $\cL_\nu$.  Let
$\Phi_\nu: X_\Si\to M_\bR$ be a moment map of $\cL_\nu$, and
let $\Psi_\nu = \Phi_\nu\circ j_0: N_\bR \to M_\bR$.
Define $\phi_\nu: M_\bR \times N_\bR\to  T^\vee_\bR \times M_\bR $
by $\phi_\nu (x,y) = (p(x), \Psi_\nu(y))$,
where $p: M_\bR \to M_\bR/M=T^\vee_\bR$ is the natural projection.
Given any line bundle $\cL_\vc$ over $X_\Si$,
$$
L_{\nu, \vc,h} :=
\phi_\nu(\bL_{\vc,h})
=\{ (p \circ \Psi_{\vc,h}(y), \Psi_\nu(y))\mid y\in N_\bR \} \subset T_\bR^\vee\times M_\bR
$$
is (in the equivalence class of) Abouzaid's Lagrangian brane
associated to the line bundle $\cL_\vc$.

Abouzaid defined a relative Fukaya category $Fuk((\bC^*)^n, M)$, where
$M$ is a fiber of $W:(\bC^*)\to \bC$, and proved that $\cL_\vc\mapsto L_{\nu,\vc,h}$
defines a full embedding
$$
DCoh(X_\Sigma) \to D^\pi Fuk((\bC^*)^n, M),
$$
which is expected to be an equivalence
when $X_\Sigma$ is Fano. So when $X_\Sigma$ is a smooth projective Fano toric variety,
it is natural to expect
\begin{equation}\label{eqn:conjecture}
DCoh(X_\Sigma) \cong   DFuk(T^*T_\bR^\vee,\bar{\Lambda}_\Si) \cong D^\pi Fuk((\bC^*)^n, M),
\end{equation}
where the equivalences are given by $\cL_\vc\mapsto \bar{\bL}_{\vc,h}\mapsto L_{\nu,\vc,h}$.

Abouzaid's work (as the authors understand it) is inspired in part by T-duality, but in it there is not the
emphasis (as there is here) that T-duality is the precise mechanism for mirror symmetry, nor is there
a connection to constructible sheaves.

\end{appendix}


\begin{thebibliography}{CK}

\bibitem[AB]{AB} M.F. Atiyah and R. Bott,
``The moment map and equivariant cohomology,''
Topology {\bf 23} (1984), no. 1, 1--28.


\bibitem[Ab1]{Ab1} M. Abouzaid, ``Homogeneous coordinate rings
and mirror symmetry for toric varieties,"
Geometry \& Topology {\bf 10} (2006), 1097--1156.


\bibitem[Ab2]{Ab2} M. Abouzaid, ``Morse homology, tropical geometry, and homological mirror symmetry for
toric varieties,"  Selecta Math. (N.S.) {\bf 15} (2009), no. 2, 189--270.

\bibitem[AP]{AP} D. Arinkin and A. Polishchuk,
``Fukaya category and Fourier transform,''
in {\em Winter School on Mirror Symmetry, Vector Bundles and Lagrangian Submanifolds,}
AMS/IP Stud. Adv. Math. {\bf 23} (2001) 261--274.

\bibitem[Au]{Au} D. Auroux, ``Mirror symmetry and T-duality in the complement of
the anticanonical divisor,"  J. G\"{o}kova Geom. Topol. GGT {\bf 1} (2007), 51--91.

\bibitem[AKO1]{AKO1} D. Auroux, L. Katzarkov, D. Orlov,
``Mirror symmetry for weighted projective planes and
their noncommutative deformations'',
Ann. of Math. (2)  {\bf 167}  (2008),  no. 3, 867--943.

\bibitem[AKO2]{AKO2} D. Auroux, L. Katzarkov, D. Orlov,
``Mirror symmetry for del Pezzo surfaces: vanishing cycles
and coherent sheaves'',  Invent. Math. {\bf 166}  (2006),  no. 3, 537--582.
%
%
\bibitem[Bo]{Bondal} A. Bondal, ``Derived categories of toric varieties,'' in
{\em Convex and Algebraic geometry, Oberwolfach conference
reports, EMS Publishing House} {\bf 3} (2006) 284--286.


\bibitem[BR]{BR} A. Bondal, W.-D. Ruan, ``Mirror symmetry for
weighted projective spaces," in preparation.
%

%

\bibitem[CL]{CL} K. Chan and N.-C. Leung, 
``Mirror symmetry for toric Fano manifolds via SYZ transformations,''
Adv. Math. {\bf 223} (2010), no. 3, 797-839. 

\bibitem[CO]{CO} C.-H. Cho and Y.-G. Oh, ``Floer cohomology
and disc instantons of Lagrangian torus fibers in Fano toric
manifolds,'' Asian J. Math {\bf 10} (2006) 773--814.
%

\bibitem[Dr]{Dr} V. Drinfeld, �DG quotients of DG categories,J. Algebra {\bf 272} (2004), no. 2, 64391.

\bibitem[F]{F}B. Fang, ``Homological mirror symmetry is $T$-duality for $\mathbb{P}^n$,"
Commun. Number Theory Phys. 2 (2008), no. 4, 719--742.

\bibitem[FLTZ]{FLTZ}B. Fang, C.-C. Liu, D. Treumann and E. Zaslow,
``A Categorification of Morelli's Theorem," {\tt arXiv:1007.0053}.

\bibitem[Ful]{Fu}W. Fulton, {\em Introduction to toric varieties,}
Annals of Mathematics Studies {\bf 131},
Princeton University Press, 1993.

\bibitem[Fuk]{Fuk}K. Fukaya,
``Mirror symmetry of abelian varieties and multi-theta functions.''
J. Algebraic Geom. {\bf 11} (2002), no. 3, 393--512.

\bibitem[FOOO1]{FOOO1} K. Fukaya, Y.-G. Oh, H. Ohta, K. Ono, ``Lagrangian
Floer theory on compact toric manifolds I,'' Duke Math J. {\bf 151} (2010), no.1, 
23-175.

\bibitem[FOOO2]{FOOO2} K. Fukaya, Y.-G. Oh, H. Ohta, K. Ono, ``Lagrangian
Floer theory on compact toric manifolds II: Bulk deformations,'' {\tt arXiv:0810:5654}.
%
%

\bibitem[Gu1]{Gu} V. Guillemin,
``Kaehler structures on toric varieties,''
J. Differential Geom. {\bf 40}  (1994),  no. 2, 285--309.

\bibitem[Gu2]{Gu2} V. Guillemin, ``Moment maps and combinatorial invariants
of Hamiltonian $T^n$-spaces,''
Progress in Mathematics, {\bf 122},
Birkh\"{a}user Boston, Inc., Boston, MA, 1994.
%

\bibitem[HV]{HV} K. Hori and C. Vafa, ``Mirror Symmetry,'' {\tt hep-th/0002222.}

\bibitem[KT]{KT} Y. Karshon and S. Tolman,
``The moment map and line bundles over pre-symplectic toric
manifolds,"  J. Differential Geom.  {\bf 38}  (1993),  no. 3, 465--484.

\bibitem[KS]{KS} M. Kashiwara and P. Schapira, ``Sheaves on Manifolds,''
Grundlehren der Mathematischen Wissenschafte {\bf 292}, Springer-Verlag, 1994.
%

\bibitem[K1]{K1} M. Kontsevich, ``Homological algebra of mirror symmetry",
Proceedings of the International Congress of Mathematicians (Z\"{u}rich, 1994),
1995, 120--139.

\bibitem[K2]{K2} M. Kontsevich, course at ENS, 1998,
{\tt http://www.math.uchicago.edu/$\sim$mitya/langlands/kontsevich.ps}

\bibitem[LYZ]{LYZ} N.-C. Leung, S.-T. Yau and E. Zaslow,
``From special Lagrangian to Hermitian-Yang-Mills via Fourier-Mukai transform,"
Adv. Theor. Math. Phys. {\bf 4} (2000) 1319--1341.


\bibitem [Mo]{M} R. Morelli,
``The K theory of a toric variety'',
Adv. Math.  {\bf 100}  (1993),  no. 2, 154--182.

\bibitem [N1]{N} D. Nadler,  ``Microlocal Branes are Constructible Sheaves,"
 Selecta Math. (N.S.) {\bf 15} (2009), no. 4, 563--619.

\bibitem [N2]{N2} D. Nadler, ``Springer theory via the Hitchin fibration,"
{\tt arXiv:0806.4566.}

\bibitem [NZ]{NZ} D. Nadler and E. Zaslow, ``Constructible sheaves and the Fukaya category," J. Amer. Math. Soc. {\bf 22}  (2009) 233--286.

\bibitem [Od]{Od} T. Oda,
{\em Convex bodies and algebraic geometry. An introduction to the theory of toric varieties.}
Translated from the Japanese. Ergebnisse der Mathematik und ihrer Grenzgebiete 3.
{\bf 15}. Springer-Verlag, Berlin, 1988.




%

\bibitem[S1]{S2} P. Seidel, ``More about vanishing cycles and mutation",
Symplectic Geometry and Mirror Symmetry, Proc. 4th KIAS International
Conference (Seoul, 2000), 429-465.

\bibitem[S2]{S3} P. Seidel, ``Fukaya categories and Picard-Lefschetz theory",
Z\"{u}rich Lectures in Advanced Mathematics.
European Mathematical Society (EMS), 2008.

\bibitem[SYZ]{SYZ} A. Strominger, S.-T. Yau, and E. Zaslow,
``Mirror symmetry is $T$-duality," Nuclear Phys. {\bf B479} (1996) 243--259.

\bibitem[Su]{Su} A. Subotic, ``Monoidal structure on the Fukaya category,'' thesis
in preparation.


\bibitem[Th]{Th} R.W. Thomason,
``The classification of triangulated subcategories,''
 Compositio Math.  {\bf 105}  (1997),  no. 1, 1--27.

\bibitem[Tr]{T} D. Treumann, ``Remarks on the nonequivariant coherent-constructible
correspondence,'' preprint.



\bibitem[U]{U} K. Ueda, ``Homological mirror symmetry for toric del Pezzo surfaces,"
Comm. Math. Phys.  {\bf 264}  (2006),  no. 1, 71--85.

\bibitem[UY]{UY} K. Ueda and M. Yamazaki, ``Homological mirror symmetry for toric orbifolds of toric del Pezzo
surfaces", {\tt arXiv:math/0703267}.

\bibitem[vdDM]{vM} L. van den Dries and C. Miller,
``Geometric categories and o-minimal structures,"
Duke Math. J. {\bf 84} (1996), no. 2, 497--540.
\end{thebibliography}
\end{document}